\documentclass{article}
\usepackage[utf8]{inputenc}
\usepackage{mathrsfs}
\usepackage{amssymb}
\usepackage{amsmath}
\usepackage{bm}
\usepackage{cite}
\usepackage{authblk}
\usepackage[algo2e]{algorithm2e}
\usepackage{tikz,pgfplots}
\pgfplotsset{compat=1.18}
\usepackage{amsfonts}   
\usepackage{subfigure}
\usepackage{graphicx}       
\usepackage{braket}
\usepackage{booktabs}
\usepackage{multirow}
\usepackage{threeparttable}
\usepackage[margin=1in]{geometry}

\usepackage[normalem]{ulem}
\usepackage[colorlinks,linkcolor=blue,anchorcolor=blue,citecolor=red]{hyperref}
\usepackage{enumerate}
\usepackage{cases}
\usepackage{mathrsfs}
\usepackage{calligra}
\newtheorem{remark}{Remark}

\newcommand\bx{\boldsymbol{x}}
\newcommand\bn{\boldsymbol{n}}
\newcommand\bu{\boldsymbol{u}}
\newcommand\bv{{\boldsymbol{v}}}

\newcommand\bc{\boldsymbol{c}}

\newcommand\bh{\boldsymbol{h}}
\newcommand\bbR{\mathbb{R}}

\newcommand\Kn{{\rm Kn}}

\newcommand\mQ{\mathcal{Q}}

\newcommand\pd[2]{\frac{\partial {#1}}{\partial {#2}}}

\newcommand\tr{\text{snap}}

\newcommand{\dd}{\; \mathrm{d}}

\graphicspath{{media/}}     

\numberwithin{equation}{section}

\title{A data-based conservative reduced collision model for the Boltzmann equation}
\date{}
\author{
  ~~Pei Zhang\footnote{Beijing Computational Science Research Center,
  Beijing, China, email:
  \texttt{zhangpei@csrc.ac.cn}.},
  ~~Yanli Wang\footnote{Beijing Computational Science Research Center,
  Beijing, China, email:
  \texttt{ylwang@csrc.ac.cn}.}
   }

\begin{document}
\maketitle

\begin{abstract}
The evaluation of the quadratic collision operator is one of the major computational bottlenecks in solving the Boltzmann equation owing to the high dimensionality and the nonlinear structure of the collision term. In this work, a conservative reduced collision (RC) model is developed within the discrete velocity framework. The bilinear collision operator is approximated through a kernel tensor preassembled with the problem-dependent reduced bases for the distribution function and the collision term. Consequently, online evaluation of the quadratic collision operator is reduced to tensor contractions. Exact discrete conservation of mass, momentum, and energy is enforced by constructing the collision basis in the orthogonal complement of the space spanned by the discrete collision invariants. Several numerical experiments validate the efficiency of this RC model. 

\end{abstract}


\section{Introduction}

Rarefied gas flows occur in a broad range of applications, including aerospace flight, vacuum technologies, and micro-/nano-scale flows, where the molecular mean free path is no longer negligible relative to the characteristic length scale \cite{shen2006rarefied}. In such regimes, continuum models based on the Navier--Stokes equations lose accuracy and eventually become invalid. As the fundamental kinetic equation for non-equilibrium gas transport, the Boltzmann equation provides a unified description across a wide range of Knudsen numbers. However, numerically solving Boltzmann equation directly is computationally expensive. The distribution function is defined on a high-dimensional phase space, while the collision term is a nonlinear quadratic integral over the velocity variables. In particular, the repeated evaluation of the collision operator constitutes a major computational bottleneck. A common way to reduce this cost is to replace the collision operator with a simplified model such as the BGK \cite{bhatnagar1954model} or the Shakhov \cite{shakhov1968generalization} model, which is much cheaper to evaluate but no longer reproduces the full collision dynamics.

Numerical methods for the Boltzmann equation can be broadly divided into stochastic and deterministic methods. Among stochastic methods, the direct simulation Monte Carlo method is widely utilized for highly non-equilibrium rarefied flows, but statistical sampling noise becomes particularly pronounced in low-speed and weakly perturbed regimes \cite{bird1994molecular, oran1998direct}. Deterministic methods, by contrast, avoid such noise and directly approximate the distribution function. A classical deterministic method is the discrete velocity method (DVM), which discretizes the velocity space using a finite set of grid points and replaces the collision integral by discrete interaction sums \cite{broadwell1964study, buet1996discrete,goldstein1989investigations, bobylev1995approximation}. Although conceptually straightforward, direct evaluation of these sums becomes prohibitively expensive as the velocity resolution increases. Fast summation techniques reduce this cost, but collision evaluation remains expensive in multidimensional and time-dependent computations \cite{mouhot2013convolutive}.

Spectral approximations constitute another important class of deterministic methods. Existing methods include Fourier spectral methods \cite{pareschi1996fourier, pareschi2000numerical, gamba2017fast}, Hermite spectral methods \cite{wang2019approximation, li2022hermite, li2023hermite,hu2020numerical}, Burnett spectral methods \cite{hu2020burnett, cai2020burnett}, polynomial Petrov–Galerkin methods \cite{kitzler2019polynomial, gamba2018galerkin}, and mapped Chebyshev spectral methods \cite{hu2022fast, hu2020petrov}. Among these methods, Fourier spectral methods are widely used because they provide high-order accuracy and express the collision operator in a convolution form. Early direct Fourier spectral methods achieved spectral accuracy but remained expensive to evaluate \cite{pareschi1996fourier, pareschi2000numerical}. Subsequent developments based on Carleman-type representations and kernel decompositions enabled fast spectral algorithms using the fast Fourier transform \cite{mouhot2006fast, wu2013deterministic}, and these algorithms were later generalized to broader classes of collision kernels \cite{gamba2017fast}. Fast spectral techniques have also been extended to linearized collision operators \cite{yin2025fast}. Nevertheless, the collision operator must be evaluated repeatedly at every spatial location and at each iteration or time step. Its evaluation can therefore still dominate the computational cost in multidimensional and unsteady simulations, as well as in steady-state simulations that require many iterations.

Recent developments in data-driven modeling have motivated the use of neural-network surrogates for the Boltzmann collision operator \cite{holloway2021acceleration, xiao2021using, miller2022neural, xiao2023relaxnet} and reduced-order collision models \cite{alekseenko2022fast}. In particular, a neural sparse representation has been developed in \cite{li2024solving}, where the reduced bases for the distribution function and the collision term are obtained from the numerical solution of the BGK equation computed by a neural network, and the conservation is enforced by projecting the collision basis. Structure-preserving neural-network models have also been developed for the closure of Boltzmann moment systems \cite{schotthofer2021structurepreserving}. These methods demonstrate the potential of such approximations to reduce computational cost while retaining important physical properties. Another closely related line of research exploits tensor decompositions and low-rank structure to reduce the cost of high-dimensional kinetic equations. Tensor methods have been developed for high-dimensional linear partial differential equations and the Boltzmann--BGK equation \cite{boelens2018parallel, boelens2020tensor}. Dynamical low-rank approximation provides a general framework for evolving time-dependent low-rank factors \cite{koch2007dynamicala, lubich2014projectorsplitting} and has been adapted to the nonlinear Boltzmann equation \cite{hu2021adaptive}. Reduced-basis methods have also been developed for related kinetic transport equations, such as the radiative transfer equation \cite{peng2022reduced}. These methods primarily reduce phase-space complexity by representing and evolving the distribution function in a compressed form. For the nonlinear Boltzmann equation, the collision term remains part of the online low-rank evolution.

In this work, we develop a conservative reduced collision (RC) model within the discrete velocity framework to reduce the cost of each evaluation of the collision operator. In the offline stage, separate bases for the distribution function and the collision term are constructed from problem-dependent offline data, and a reduced kernel tensor is assembled to approximate the bilinear collision operator with the fast Fourier spectral method (FSM). Online collision evaluation then reduces to quadratic tensor contractions with complexity $\mathcal{O}(N_f^2N_Q)$, where $N_f$ and $N_Q$ denote the rank numbers of the distribution function and the collision model, respectively. To enforce the discrete conservation structure, the collision basis is projected onto the orthogonal complement of the discrete collision-invariant space, so that mass, momentum, and energy are preserved exactly at the discrete level.

Since that only the collision evaluation is modified, the model can be incorporated into existing deterministic Boltzmann solvers. The numerical assessment covers the Bobylev--Krook--Wu (BKW) relaxation problem, steady Fourier and Couette flow problems, two-dimensional lid-driven cavity flow, and an unsteady fluid diffusion problem. With the reduced ranks selected from the offline data, the proposed model remains in close agreement with the reference FSM solutions while substantially reducing the computational cost of collision evaluation.

The rest of the paper is organized as follows. The Boltzmann equation and the fast Fourier spectral method are reviewed in Sec. \ref{sec:Bolt}. The construction of the reduced collision model is presented in Sec. \ref{sec:method}, including building the conservative basis function, assembly of the offline tensor, and the procedure of the rank selection. The numerical results are shown in Sec. \ref{sec:num}, with the concluding remarks listed in Sec. \ref{sec:conclusion} and several appendices in App. \ref{app:bilinear}.

\section{Boltzmann equation and discretization of the microscopic velocity space}
\label{sec:Bolt}
In this section, the Boltzmann equation is briefly reviewed, and the discretization of the Boltzmann equation in the microscopic velocity space is introduced, where the discrete velocity method with the fast Fourier spectral discretization is utilized to approximate the quadratic collision model. 

\subsection{The Boltzmann equation and collision kernel}
The dynamics of a rarefied gas are governed by the Boltzmann equation, which describes the time evolution of the particle distribution function $f(t, \bx, \bv)$ in phase space. Let $t \geqslant 0$ denote time, $\bx \in \Omega \subset \mathbb{R}^3$ the spatial coordinate, and $\bv \in \mathbb{R}^3$ the microscopic velocity. The equation reads:
\begin{equation}
    \label{eq:Boltzmann}
    \pd{f}{t} + \bv \cdot \nabla_{\bx} f  = \mQ(f, f),  \qquad t \geqslant 0, \quad \bx \in \Omega, \quad \bv \in \bbR^{3},
\end{equation}
where the collision operator $\mathcal{Q}(f, f)$ models the binary interactions between particles and has the following integral form:
\begin{equation}
    \label{eq:quic_coll}
    \mathcal{Q}(f, f)=\int_{\bbR^{3}} \int_{\mathbb{S}^{2}} B\left(\bv-\bv_{{\ast}}, \sigma \right)\left[f\left(\bv_{{\ast}}^{\prime}\right) f\left(\bv^{\prime}\right)-f\left(\bv_{{\ast}}\right) f(\bv)\right] \dd \sigma \dd \boldsymbol{v}_{{\ast}},
\end{equation}
where $\bv$ and $\bv_{\ast}$ denote the pre-collision velocities of two particles, while $\bv^{\prime}$ and $\bv_{\ast}^{\prime}$ represent their post-collision velocities, which are determined by the conservation laws of momentum and energy, 
\begin{equation}
    \bv' = \frac{\bv+\bv_{\ast}}{2}+\frac{|\bv-\bv_{\ast}|}{2}\sigma,\qquad \bv'_{\ast} = \frac{\bv+\bv_{\ast}}{2} -\frac{|\bv-\bv_{\ast}|}{2}\sigma,
\end{equation}
where $\sigma\in\mathbb{S}^2$ is a unit vector.
The collision kernel $B(\cdot,\cdot )$ is a non-negative function which is determined by the intermolecular force between gas molecules. In this work, the variable hard-sphere model is adopted \cite{bird1994molecular}, where the collision kernel $B$ is independent of $\sigma$ as 
\begin{equation}
\label{eq:B}
    B = C_{\alpha} |\bv - \bv_{\ast}|^{\alpha}, \quad \alpha = 2(1-\omega).
\end{equation}
Here, the parameter $\omega$ is the viscosity index of the gas, and the constant $C_\alpha$ is related to the diameter of the gas molecule. The collision operator conserves the mass, momentum, and energy as 
\begin{equation}
\label{eq:col_con}
    \int_{\mathbb{R}^3} \mQ(f,f)\,\varphi(\bv)\,\dd \bv = 0, \qquad \varphi(\bv)=1, \bv, |\bv|^2.
\end{equation}
Besides, the macroscopic observables, including density $\rho$, macroscopic velocity $\bu$, and temperature $T$, can be obtained from the distribution function as 
\begin{equation}
\label{eq:mac_var}
    \rho = \int_{\bbR^3} f \dd \bv,\quad \bu = \frac{1}{\rho}\int_{\bbR^3}\bv f\dd \bv,\quad T = \frac{1}{3\rho}\int_{\bbR^3}|\bv-\bu|^2 f\dd \bv.
\end{equation}

The efficient simulation of the Boltzmann equation is challenging due to its high dimensionality, the complex quadratic collision model, etc. Following the method in \cite{gamba2017fast}, the fast Fourier spectral method is utilized here to approximate the collision term. Precisely, the distribution function $f(\bv)$ is expanded in Fourier form as  
\begin{equation}
\label{eq:analy_Fourier}
f(\bv) = \frac{1}{(2\pi)^3} \int_{\mathbb{R}^3} \widehat{f}(\boldsymbol{\xi}) e^{i \bv \cdot \boldsymbol{\xi}} \, \dd\boldsymbol{\xi}, \quad \widehat{f}(\boldsymbol{\xi}) = \int_{\mathbb{R}^3} f(\bv) e^{-i \bv \cdot \boldsymbol{\xi}} \, \dd\bv.
\end{equation}
Here, the variables $(t, \bx)$ are omitted for simplicity. Thus, the collision term leads to a convolution-type representation in the Fourier space as 
{\small
\begin{equation}
    \begin{aligned}
        &\widehat{Q}(\boldsymbol{\xi}) = \int_{\mathbb{R}^3} \int_{\mathbb{R}^3} \int_{\mathbb{S}^2} e^{-i \bv \cdot \boldsymbol{\xi}} B(|\bv - \bv_*|) [f(\bv')f(\bv'_*) - f(\bv)f(\bv_*)]  \dd\sigma  \dd \bv_*  \dd \bv \\
        & = \frac{1}{(2\pi)^6} \int_{\mathbb{R}^3} \int_{\mathbb{R}^3} \int_{\mathbb{R}^3} 
        \int_{\mathbb{R}^3}\int_{\mathbb{S}^2} B(|\bv - \bv_*|) \widehat{f}(\boldsymbol{\eta}) \widehat{f}(\boldsymbol{\zeta})\left( e^{i \bv' \cdot \boldsymbol{\eta}} e^{i \bv'_* \cdot \boldsymbol{\zeta}} - e^{i \bv \cdot \boldsymbol{\eta}} e^{i \bv_* \cdot \boldsymbol{\zeta}} \right) e^{-i \bv \cdot \boldsymbol{\xi}}  \dd\sigma  \dd\bv_*  \dd\bv \dd \boldsymbol{\eta} \dd\boldsymbol{\zeta}.
    \end{aligned}
\end{equation}
}

\subsection{Discrete velocity method and the fast Fourier spectral approximation}
\label{sec:dvm}
For the discrete velocity method (DVM), the microscopic velocity domain is truncated to the bounded region $\Omega_v = [-L, L]^3$, which is large enough to capture the essential dynamics of the molecular velocity distribution. Assuming that $f$ is supported in the ball $\mathcal{B}_R$ of radius $R$, the relative velocity $\bv-\bv_*$ is truncated to $\mathcal{B}_{2R}$, and $L\geqslant\frac{3+\sqrt{2}}{2}R$ is chosen to avoid aliasing \cite{pareschi2000numerical}. To apply the discrete fast Fourier spectral method, the discrete velocity set is chosen as 
\begin{equation}
\label{eq:dis_v}
    \bv_{\boldsymbol{p}}=\frac{L}{N}\boldsymbol{p},
    \qquad
    \Delta v=\frac{L}{N},
\end{equation}
where $\bm{p} \in \mathcal{I}_N = \{-N,\ldots, N-1\}^3$ is a multi-index with $N$ the number of Fourier mode. Thus, the discrete distribution function is approximated in the framework of the Fourier spectral method as 
\begin{equation}
\label{eq:Fourier_dis_f}
    f_N(\bv) = \frac{1}{(2L)^3}\sum_{\boldsymbol{k}\in\mathcal{I}_N}  \widehat f_{\boldsymbol{k}}
    \exp\left(\frac{\mathrm{i}\pi}{L}\boldsymbol{k}\cdot\bv \right),
\end{equation}
where $\boldsymbol{k}$ is a Fourier frequency index, and the Fourier coefficients $\widehat f_{\boldsymbol{k}}$ are calculated as 
\begin{equation}
\label{eq:coe_f}
    \widehat f_{\boldsymbol{k}} =(\Delta v)^3\sum_{\boldsymbol{p}\in\mathcal{I}_N}  f_{\boldsymbol{p}}   \exp\left(
        -\frac{\mathrm{i}\pi}{N}  \boldsymbol{k}\cdot\boldsymbol{p}  \right),
    \qquad
    \boldsymbol{k}\in\mathcal{I}_N,
\end{equation}
where
$f_{\boldsymbol{p}}=f(\bv_{\boldsymbol{p}})$ and the variables $(t,\bx)$ are omitted for simplicity. Thus, the Fourier coefficients of the collision term are evaluated as
\begin{equation}
\label{eq:dis_hatQ}
\begin{aligned}
    \widehat{Q}_{\boldsymbol{k}} &=
    \frac{2\pi^2}{L^3}\sum_{n=1}^{J}\sum_{m=1}^{M}w_n^{\varrho} w_m^{\sigma} B(\varrho_n)
    \operatorname{sinc}\left(\frac{\pi \varrho_n |\boldsymbol{k}|}{2L}\right)I_{n,\boldsymbol{k},m}   \\
    &\qquad \qquad  -\frac{2\pi^2}{L^3}\sum_{\boldsymbol{\ell}\in\mathcal{I}_N}
    \widehat f_{\boldsymbol{k}-\boldsymbol{\ell}}  \widehat f_{\boldsymbol{\ell}}
    \left(\sum_{n=1}^{J} w_n^{\varrho} B(\varrho_n) \operatorname{sinc} \left(\frac{\pi \varrho_n |\boldsymbol{\ell}|}{L} \right) \right),
    \qquad    \boldsymbol{k}\in\mathcal{I}_N,
\end{aligned}
\end{equation}
where
\begin{equation}
\label{eq:coe_L}
\begin{aligned}
    I_{n,\boldsymbol{k},m}  = \sum_{\boldsymbol{\ell}\in\mathcal{I}_N}
    \widehat f_{\boldsymbol{k}-\boldsymbol{\ell}}\widehat f_{\boldsymbol{\ell}}
    \bigg[ & \cos\left(\frac{\pi \varrho_n\boldsymbol{\sigma}_m
    \cdot (\boldsymbol{k}-\boldsymbol{\ell})}{2L} \right)
    \cos\left( \frac{\pi \varrho_n  \boldsymbol{\sigma}_m\cdot \boldsymbol{\ell}}{2L}  \right)  \\
    &+ \sin\left( \frac{\pi \varrho_n \boldsymbol{\sigma}_m\cdot
        (\boldsymbol{k}-\boldsymbol{\ell})}{2L} \right)
    \sin\left(  \frac{\pi \varrho_n \boldsymbol{\sigma}_m\cdot  \boldsymbol{\ell}}{2L}  \right) \bigg].
\end{aligned}
\end{equation}
Here, $\{\varrho_n, w_n^{\varrho}\}_{n=1}^{J}$ are the nodes and weights of the Jacobi--Gauss--Radau quadrature on $[0,2R]$, and $\{\boldsymbol{\sigma}_m, w_m^{\sigma}\}_{m=1}^{M}$ are the nodes and weights of the Lebedev quadrature on $\mathbb{S}^2$. Besides, $B(\varrho_n)$ denotes the collision kernel \eqref{eq:B} evaluated at $|\bv-\bv_*|=\varrho_n$, and $\operatorname{sinc}(x)=\sin x/x$.

Moreover, the discrete collision terms utilized in the DVM method are obtained by transforming the Fourier coefficients $\widehat Q_{\boldsymbol{k}}$ back to the phase space 
\begin{equation}
\label{eq:dis_Q}
    Q_{\boldsymbol{p}} = \frac{1}{(2L)^3}\sum_{\boldsymbol{k}\in\mathcal{I}_N}
    \widehat Q_{\boldsymbol{k}}\exp\left(\frac{\mathrm{i}\pi}{N}\boldsymbol{k}\cdot\boldsymbol{p}\right),
    \qquad \boldsymbol{p}\in\mathcal{I}_N.
\end{equation}
With the discrete velocity set \eqref{eq:dis_v}, the total number of velocity points is $N_v = (2N)^3$. To build the governing equations of the Boltzmann equation in the framework of DVM, the multi-index $\bm{p}$ is first flattened into a single index $s, 1\leqslant s \leqslant N_v$. Then, the Boltzmann equation \eqref{eq:Boltzmann} is reduced into a coupled system of partial differential equations as
\begin{equation}
\label{eq:dis_Bol}
\frac{\partial f_s}{\partial t} + \bv_s \cdot \nabla_{\bx} f_s = Q_s(\boldsymbol{f},\boldsymbol{f}),
\quad s = 1, 2, \dots, N_v,
\end{equation}
where $\bv_s=(v_{1,s},v_{2,s},v_{3,s})^T$ denotes the discrete velocity point with index $s$, and $Q_s(\boldsymbol{f},\boldsymbol{f})$ denotes the discrete collision term at $\bv_s$, which depends on the discrete distribution function $\boldsymbol{f}=(f_1,\dots,f_{N_v})^T\in\mathbb{R}^{N_v}$. The discrete collision term is defined as
\begin{equation}
    \boldsymbol{Q}=\big(Q_1(\boldsymbol{f},\boldsymbol{f}),\dots,Q_{N_v}(\boldsymbol{f},\boldsymbol{f})\big)^T\in\mathbb{R}^{N_v}.
\end{equation}
The system \eqref{eq:dis_Bol} is then utilized as the semi-discrete system to numerically solve the Boltzmann equation. The fast Fourier spectral method is a well-established method for the Boltzmann equation, and we refer \cite{gamba2017fast} for more details. Here, we want to emphasize that the 
resulting fast spectral evaluation \eqref{eq:dis_hatQ} has computational complexity $\mathcal{O}(MN^4\log N)$ \cite{gamba2017fast}. Based on this fast Fourier spectral method, a reduced formulation for the collision term is developed in the next section.

\section{Constructing the conservative reduced collision model}
\label{sec:method}

In this section, a conservative reduced model for the collision is developed, based on the discrete FSM presented in Sec. \ref{sec:Bolt}, where the complex computations of the collision operator in the FSM are replaced by a reduced bilinear form, with compatibility with the discrete conservation structure of the Boltzmann equation retained.

\subsection{Reduced collision model}
\label{subsec:rc-ansatz}
With the DVM discretization introduced in Sec. \ref{sec:Bolt}, the discrete distribution function $\boldsymbol{f}$ and the discrete collision term $\boldsymbol{Q}$ are both vectors in $\mathbb{R}^{N_v}$. Assuming that we have already found a series of basis functions, which can approximate the distribution function $\bm{f}$ and the quadratic collision term $\bm{Q}$ with high efficiency, the reduced collision model is built following the thought in the spectral method. Precisely, let
\begin{equation}
    \label{eq:basic_fun}
    \{\boldsymbol{g}_r\}_{r=1}^{N_f}\subset\mathbb{R}^{N_v}, \qquad \{\boldsymbol{h}_r\}_{r=1}^{N_Q}\subset\mathbb{R}^{N_v}
\end{equation}
be the orthogonal bases for the discrete distribution functions and the discrete collision term, respectively. Then $\bm{f}$ and $\bm{Q}$ can be approximated as 
\begin{equation}
\label{eq:rc-expansion}
  \boldsymbol{f} \approx \sum_{r=1}^{N_f} \hat f_r\, \boldsymbol{g}_r, \qquad
  \boldsymbol{Q} \approx \sum_{r=1}^{N_Q} \hat q_r\, \boldsymbol{h}_r,
\end{equation}
where the corresponding expansion coefficients are obtained by the orthogonal projection as 
\begin{equation}
    \label{eq:coeff}
    \hat{f}_r = \left\langle \boldsymbol{f},\boldsymbol{g}_r \right \rangle,\qquad \hat{q}_r = \left\langle \boldsymbol{Q}, \boldsymbol{h}_r \right \rangle,
\end{equation}
with inner product $\langle \bm{a},\bm{b}\rangle=\sum_{s=1}^{N_v} a_s b_s$ for $\bm{a},\bm{b}\in\bbR^{N_v}$. Substituting the reduced expansion \eqref{eq:rc-expansion} of $\boldsymbol{f}$ into $\mQ(\boldsymbol{f},\boldsymbol{f})$ yields
\begin{equation}
\label{eq:coe_Q}
    \hat{q}_r = \sum_{i=1}^{N_f}\sum_{j=1}^{N_f} \hat{f}_i \hat{f}_j \left\langle \mQ(\boldsymbol{g}_i,\boldsymbol{g}_j),\boldsymbol{h}_r \right \rangle,
\end{equation}
where $\mQ(\boldsymbol{g}_i,\boldsymbol{g}_j)$ is obtained by FSM and the detailed results are listed in App. \ref{app:bilinear}. Defining the reduced kernel tensor $\mathcal{K}$ as 
\begin{equation}
\label{eq:kerner-tensor}
    \mathcal{K}_{ijr} = \big\langle \mQ(\boldsymbol{g}_i,\boldsymbol{g}_j),\, \boldsymbol{h}_r\big\rangle\in \bbR^{N_f\times N_f\times N_Q},
\end{equation}
for any discrete distribution function, the reduced collision model is then built as 
\begin{equation}
\label{eq:redoce_col}
    \mQ_{\text{RC}}(\boldsymbol{f},\boldsymbol{f}) = \sum_{r=1}^{N_Q} \left(\sum_{i=1}^{N_f}\sum_{j=1}^{N_f} \mathcal{K}_{ijr} \hat{f}_i \hat{f}_j \right)\boldsymbol{h}_r.
\end{equation}
Since the reduced kernel tensor $\mathcal{K}$ depends only on the reduced bases $\bm{g}_r$ and $\bm{h}_r$, it is independent of the discrete distribution function $\boldsymbol{f}$ and can be precomputed offline. 
With the well-prepared reduced kernel tensor $\mathcal{K}$, the computational cost of evaluating the collision model online reduces to $\mathcal{O}(N_f^2N_Q)$, which is greatly reduced compared to FSM. 

Here, the only point left is how to find the highly efficient basis functions $\bm{g}_r$ and $\bh_r$, which will be introduced in the following section. 



\subsection{Construction of the conservative reduced basis}
\label{subsec:offline}
The precision and physical consistency of the reduced collision model depend entirely on the reduced basis $\boldsymbol{g}_r$ and $\boldsymbol{h}_r$. Instead of choosing the classical basis functions in the spectral method, a series of basis vectors are built on the known numerical solution in this work.

\paragraph{Reduced-space construction from snapshot matrices}
We do not aim to construct a universal basis for the Boltzmann equation, but propose a general numerical method to build the problem-dependent basis. Assuming that a numerical solution with lower resolution is already obtained for the target problem, a series of discrete distribution functions at $(t, \bx)_z$ are selected from this obtained numerical solution as 
\begin{equation}
    \label{eq:tr_f}
    \boldsymbol{f}^{\tr}_z \triangleq \boldsymbol{f}^{\tr}(t_z, \bx_z),  \qquad 
    \qquad z=1,\dots,N_z. 
\end{equation}
Then, the corresponding discrete collision terms are obtained as 
\begin{equation}
    \label{eq:tr_Q}
    \boldsymbol{Q}^{\tr}_z = \mQ(\boldsymbol{f}^{\tr}_z,\boldsymbol{f}^{\tr}_z), \qquad z = 1, \cdots, N_z,
\end{equation}
with the corresponding snapshot matrices $\mathbf{F}^{\tr}$ and $\mathbf{Q}^{\tr}$ written as
\begin{equation}
\label{eq:snap_matrix}
\begin{aligned}
    \mathbf{F}^{\tr} = [\boldsymbol{f}_1^{\tr},\cdots,\boldsymbol{f}^{\tr}_{N_{z}}] \in \bbR^{N_v \times N_z},\qquad 
    \mathbf{Q}^{\tr} = [\boldsymbol{Q}_1^{\tr},\cdots,\boldsymbol{Q}^{\tr}_{N_{z}}] \in \bbR^{N_v \times N_z}.
\end{aligned}
\end{equation}
Here, the lower-resolution numerical solution can be obtained on a coarser spatial mesh, or with a larger time step, or with a simplified collision model such as the BGK or Shakhov collision model. No matter with which method this numerical solution is obtained, the computational cost should be much smaller than that of the aimed numerical solution.

To obtain the data-based basis vector, the snapshot matrices $\mathbf{F}^{\tr}$ and $\mathbf{Q}^{\tr}$ are utilized and decomposed by singular value decomposition (SVD):
\begin{equation}
\label{eq:svd}
  \mathbf{F}^{\tr}= \mathbf{U}_f \mathbf{S}_f \mathbf{V}_f^T,\qquad
  \mathbf{Q}^{\tr}= \mathbf{U}_Q \mathbf{S}_Q \mathbf{V}_Q^T,
\end{equation}
where $\mathbf{U}_f,\mathbf{U}_Q\in\mathbb{R}^{N_v\times N_z}$ are the left singular vector matrices, $\mathbf{S}_f,\mathbf{S}_Q\in\mathbb{R}^{N_z\times N_z}$ are diagonal matrices of the corresponding singular values arranged in non-increasing order, and $\mathbf{V}_f,\mathbf{V}_Q\in\mathbb{R}^{N_z\times N_z}$ are the right singular vector matrices. Here, the left singular vector matrices $\mathbf{U}_f,\mathbf{U}_Q\in\mathbb{R}^{N_v\times N_z}$ are utilized to construct the basis vectors. The basis vectors $\boldsymbol{g}_r, r = 1, \cdots, N_f$ adopted in \eqref{eq:basic_fun} are identified with the $r$-th leading columns of $\mathbf{U}_f$, while the collision basis $\bm{h}_r, r = 1, \cdots, N_Q$ are taken as the $r$-th leading columns of $\mathbf{U}_Q$. Since the columns of $\mathbf{U}_f$ and $\mathbf{U}_Q$ are orthonormal by construction, no further orthonormalization is needed. Therefore, the basis vectors can be written in matrix form as 
\begin{equation}
    \label{eq:basis_matrix}
  \mathbf{G}=[\boldsymbol{g}_1,\dots,\boldsymbol{g}_{N_f}], \qquad
  \mathbf{H}=[\boldsymbol{h}_1,\dots,\boldsymbol{h}_{N_Q}].
\end{equation}
Here, $\mathbf{G}$ is adopted directly as the reduced basis for the discrete distribution function, while the basis constructed by $\mathbf{H}$ for the discrete collision term fails to satisfy the discrete conservation constraints, which will be enforced on $\mathbf{H}$ in the following part.

\paragraph{Conservation enforcement via null-space projection} Due to the conservative property of the collision model, it is expected that the discrete collision basis vectors can also preserve the conservation of the density, momentum, and energy. Moreover, a matrix $\widetilde{\mathbf{H}}$ has this property if and only if it satisfies \eqref{eq:cons}.
\begin{equation}
\label{eq:cons}
    \mathbf{M}^T\widetilde{\mathbf{H}}=0, \qquad  \mathbf{M}=\big[\mathbf{1},\,\bm{v}^{(1)},\,\bm{v}^{(2)},\,\bm{v}^{(3)},\,|\bm{V}|^2\big] \in\bbR^{N_v\times 5},
\end{equation}
where $\mathbf{M}$ is the discrete collision-invariant matrix, the columns of which correspond to the discrete invariants associated with mass, momentum, and energy, namely
\begin{equation}
   \qquad 
    \bm{v}^{(i)}=(v_{i,1},\dots,v_{i,N_v})^T\in\bbR^{N_v\times 1}, \qquad i=1,2,3, \qquad |\bm{V}|^2=(|\bv_1|^2,\dots,|\bv_{N_v}|^2)^T\in\bbR^{N_v\times 1}.
\end{equation}
To obtain a conservative collision basis, it is constructed by projecting the preliminary collision basis $\mathbf{H}$ onto the orthogonal complement of the discrete collision-invariant \cite{gamba2009spectral,pareschi2022moment}. To achieve this, we first find an orthonormal basis matrix $\overline{\mathbf{M}}\in\mathbb{R}^{N_v\times 5}$ for the space spanned by the columns of $\mathbf{M}$, since its columns are not orthonormal. Thus,  the $\overline{\mathbf{M}}$ satisfies 
\begin{equation}
\label{eq:orth_M}
    \overline{\mathbf{M}}^{T}\overline{\mathbf{M}}=\mathbf{I}.
\end{equation}
Projecting the basis matrix of the collision term $\mathbf{H}$ \eqref{eq:basis_matrix} onto the orthogonal complement of $\overline{\mathbf{M}}$, and then orthonormalizing its columns gives 
\begin{equation}
    \widetilde{\mathbf{H}}=\Pi( \mathbf{H}-\overline{\mathbf{M}}\,\overline{\mathbf{M}}^{T}\mathbf{H}).
\end{equation}
Here, the operator $\Pi( \mathbf{H})$ means orthonormalizing the columns of the matrix $\mathbf{H}$. 
Then, it is easy to verify that 
\begin{equation}
\label{eq:pro_H}
    \overline{\mathbf{M}}^{T}\widetilde{\mathbf{H}}=0.
\end{equation}
In this case, each column of $\widetilde{\mathbf{H}}$ lies in the orthogonal complement of the discrete collision-invariant space. Consequently, any reduced collision term expanded on this basis exactly satisfies the discrete conservation constraints for mass, momentum, and energy. In the following section, $\widetilde{\mathbf{H}}$ is written as $\mathbf{H}$ for short. 

For now, the reduced basis $\mathbf{G}$ and $\mathbf{H}$ for the discrete function and the discrete collision term are built, respectively. The only point left is how to decide the reduced ranks $N_f$ and $N_Q$ based on the SVD \eqref{eq:svd}.



\begin{remark}
    A similar reduced collision model, including the SVD-based bases and the conservative projection of the collision basis, has been proposed in \cite{li2024solving} within the neural sparse representation framework, where the bases are built from the BGK solution computed by a neural network. However, the numerical error there is dominated by the approximation and training errors of the neural network, and the accuracy of the reduced collision model itself is not examined separately. In this work, the reduced collision model directly replaces FSM in deterministic solvers, so that its accuracy and efficiency can be systematically assessed against FSM in Sec.~\ref{sec:num}.
\end{remark}

\subsection{Rank selection criteria}
\label{subsec:ranksel}
Selecting the ranks $(N_f, N_Q)$ of the reduced bases is essential for balancing approximation accuracy and online computational cost. Here, the cumulative energy of the singular values is a standard criterion for quantifying the truncation error of a low-rank approximation. Let $\sigma_i^{(f)}=\mathbf{S}_{f,ii}$ and $\sigma_i^{(Q)}=\mathbf{S}_{Q,ii}$ denote the singular values of $\mathbf{F}^{\tr}$ and $\mathbf{Q}^{\tr}$, arranged in non-increasing order. The residual energies are defined as
\begin{equation}
\label{eq:energy}
E_f(N_f) = 1-\frac{\sum\limits_{i=1}^{N_f}\left(\sigma_i^{(f)}\right)^2}{\|\mathbf{F}^{\tr}\|_F^2}, \qquad 
E_Q(N_Q) = 1-\frac{\sum\limits_{i=1}^{N_Q}\left(\sigma_i^{(Q)}\right)^2}{\|\mathbf{Q}^{\tr}\|_F^2},
\end{equation}
where $\|\cdot\|_F$ denotes the Frobenius norm. The preferred rank numbers $N_f$ and $N_Q$ are chosen to satisfy
\begin{equation}
E_f(N_f)\leqslant \tau_f,\qquad E_Q(N_Q)\leqslant \tau_Q,
\end{equation}
where $\tau_f$ and $\tau_Q$ are problem-dependent tolerances to select the ranks $N_f$ and $N_Q$, respectively.

Here, this selection procedure is discussed in detail for two representative cases, the spatially homogeneous BKW problem and the steady-state Fourier flow problem, which together span the two broad problem classes considered in Sec.~\ref{sec:num}. The variation of the residual energy with the increase of rank number on different numerical settings is studied. 

\begin{figure}[!hptb]
\centering
\subfigure[$E_f(N_f)$]{
            \includegraphics[width = 0.3\textwidth]{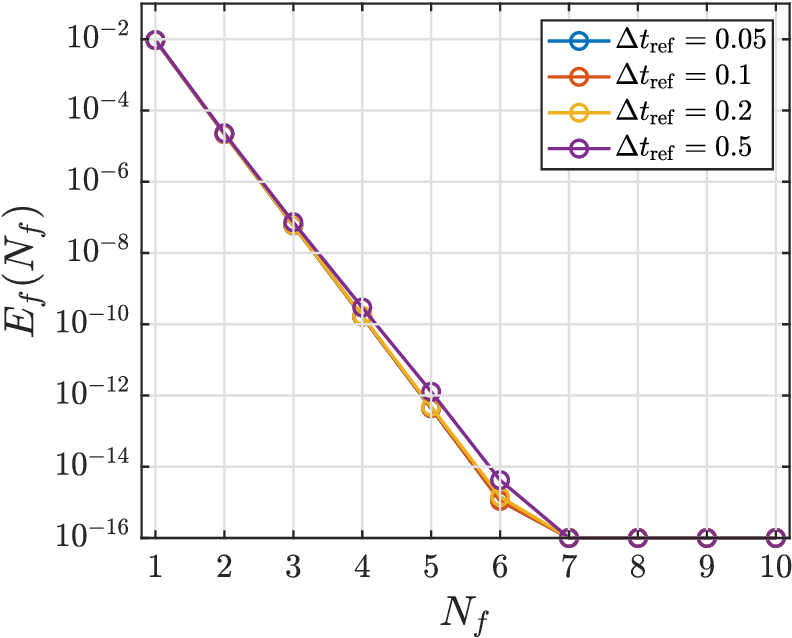}
} 
\qquad
\subfigure[$E_Q(N_Q)$]{               
            \includegraphics[width = 0.3\textwidth]{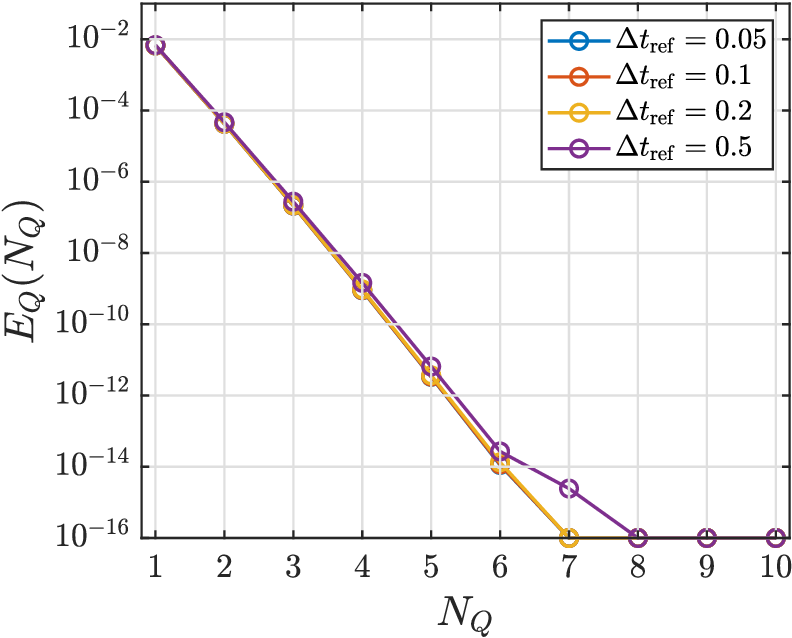}
}      
\caption{Variation of the residual energy for different rank numbers for the BKW problem with different time step length $\Delta t$. (a) $E_f(N_f)$ for the distribution function; (b) $E_Q(N_Q)$ for the collision term.}
\label{fig:bkw_energy}
\end{figure}

For the BKW problem, the final time is set to $t = 1$, and different time step lengths $\Delta t = 0.5, 0.2, 0.1$, and $0.05$ are tested. For these different sets of numerical solutions, the residual energies $E_f(N_f)$ and $E_Q(N_Q)$ with different rank numbers are shown in Fig. \ref{fig:bkw_energy}. It indicates that both energies decay rapidly and reach machine precision after only a few modes. Moreover, the behavior of both energies is similar for different time step lengths. This means that the decay of the energies only weakly depends on $\Delta t$.

\begin{figure}[!hptb]
\centering
\subfigure[$E_f(N_f)$]{
            \includegraphics[width = 0.3\textwidth]{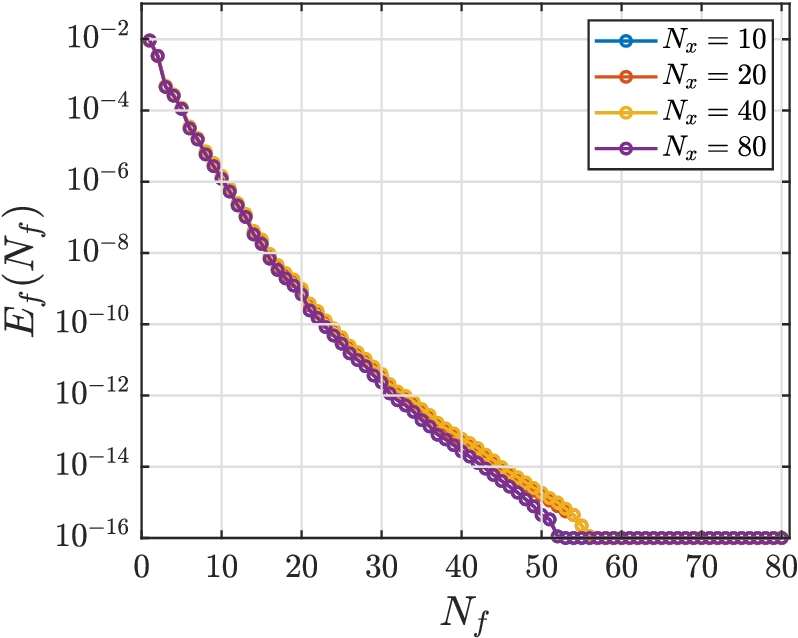}
}\qquad 
\subfigure[$E_Q(N_Q)$]{               
            \includegraphics[width = 0.3\textwidth]{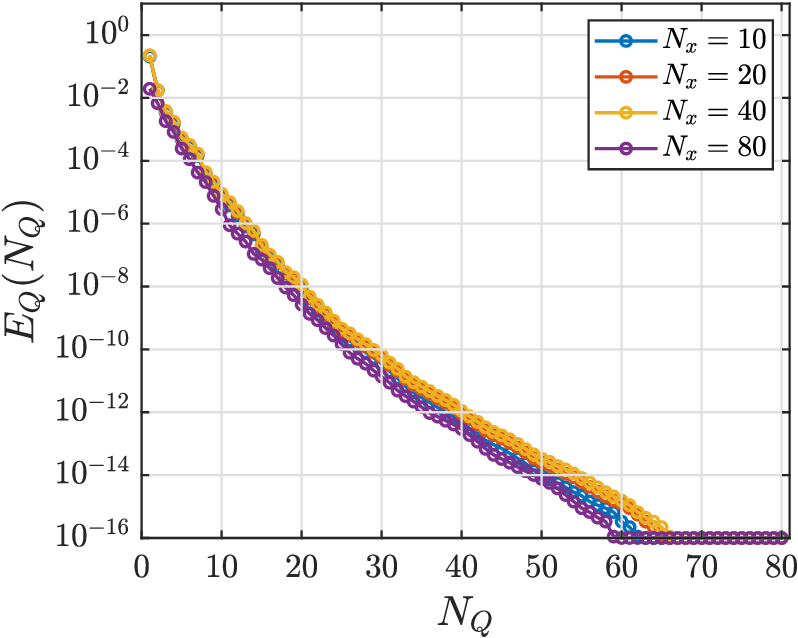}
}      
\caption{Variation of the residual energy for different rank numbers for the Fourier flow problem with different mesh sizes $N_x$. (a) $E_f(N_f)$ for the distribution function; (b) $E_Q(N_Q)$ for the collision terms.}
\label{fig:fourier_energy}
\end{figure}

For the steady-state Fourier flow problem, different numbers of spatial cells $N_x = 10, 20, 40$ and $80$ are tested.
For these different sets of numerical solutions, the residual energies $E_f(N_f)$ and $E_Q(N_Q)$ with different rank numbers are shown in Fig. \ref{fig:fourier_energy}. In this case, the decay is slower but still reaches the numerical floor, again with limited sensitivity to $N_x$.

These two examples indicate that selecting the ranks based on the residual energy in this way provides an inexpensive, qualitative guide to the required ranks. Moreover, the numerical settings such as the time step length and the mesh size play a small effect on the variation for both residual energies. The resulting choices are further validated in Sec.~\ref{sec:num}.


\section{Numerical experiments}
\label{sec:num}
In this section, the reduced collision model \eqref{eq:redoce_col} is validated with several numerical experiments, with emphasis on accuracy and computational cost relative to the FSM. We begin with the spatially homogeneous BKW solution, which provides an exact solution for evaluating the reduced collision approximation. 
The one-dimensional Fourier and Couette flow problems are considered together with the two-dimensional lid-driven cavity flow to examine the performance of the reduced collision model in spatially inhomogeneous boundary-value problems. Finally, an unsteady fluid-diffusion problem with sharp gradients is utilized to validate the new collision model in time-dependent problems. Unless otherwise stated, the velocity domain is truncated to $\Omega_v=[-L,L]^3$ with
\begin{equation}
    L=\frac{3+\sqrt{2}}{2}R,\qquad R=4.
\end{equation}
For the collision kernel \eqref{eq:B}, the constant is $C_\alpha=\frac{1}{4\pi}$, and the Maxwell molecules with $\omega=1$ are adopted.
All the simulations in Sec.~\ref{sec:BKW} and Sec.~\ref{subsec:steady} are carried out on an Intel Core i7-14650HX CPU with $16$ cores, while the time-dependent problem in Sec.~\ref{sec:diffusion} is computed on an NVIDIA A100 GPU.

\subsection{Spatially homogeneous relaxation: the BKW solution}
\label{sec:BKW}
For Maxwell molecules with collision kernel $B=\frac{1}{4\pi}$, the classical BKW problem provides an analytical solution of the spatially homogeneous Boltzmann equation
\begin{equation}
    \label{eq:bkw}
    \frac{\partial f(t, \bv)}{\partial t} = \mQ(f(t, \bv), f(t, \bv)),\qquad t>0, \bv \in \mathbb{R}^3.
\end{equation}
The analytical solution is given by~\cite{bobylev1975exact, krook1977exact}:
\begin{equation}
    f_{\rm ex}(t) = \frac{1}{2(2\pi K(t))^{3/2}} \exp \left( -\frac{\bv^2}{2K(t)} \right) \left( \frac{5K(t) - 3}{K(t)} + \frac{1 - K(t)}{K^2(t)} \bv^2 \right),
\end{equation}
where the temperature parameter evolves as 
\begin{equation}
    \label{eq:K}
K(t) = 1- \frac{2}{5}\exp(-t/6).    
\end{equation}
The corresponding exact collision term is obtained from \eqref{eq:bkw} as
\begin{equation*}
    \begin{aligned}
        \mQ_{\rm ex}(t) =  \left\{ \left( \frac{-3}{2K} + \frac{\bv^2}{2K^2} \right) f_{\rm ex} + \frac{1}{2(2\pi K)^{3/2}} \exp \left( -\frac{\bv^2}{2K} \right) \left( \frac{3}{K^2} + \frac{K - 2}{K^3} \bv^2 \right) \right\} K',\qquad K'(t) = \frac{1}{15}e^{-t/6}.
    \end{aligned}
\end{equation*}
Since the analytical solution is available, the reduced collision model is first tested on the exact distribution function $f_{\rm ex}$. In the simulation, the microscopic velocity space is discretized with $N = 24$ ($2N=48$ points per direction). The collision operator is evaluated by the FSM, with $M = 14$ Lebedev points for the angular quadrature and $J = 2N= 48$ Jacobi--Gauss--Radau points for the radial quadrature in \eqref{eq:dis_hatQ}. The offline data in \eqref{eq:tr_f} are the distribution functions obtained by solving \eqref{eq:bkw} with the fourth-order Runge--Kutta method and the reference time step length $\Delta t_{\rm ref}$ up to $t = 6$. The offline data in \eqref{eq:tr_Q} are the collision terms of these distribution functions, which are computed by the FSM. Let $\mQ_{\mathrm{FSM}}$ denote the collision term computed by the FSM, and $\mQ_{\mathrm{RC}}^{(\Delta t_{\rm ref},r)}$ denote the reduced collision term built from these basis data with the ranks $N_f=N_Q=r$. Both of them are evaluated on the exact solution $f_{\rm ex}$, and the collision errors of RC and of the FSM are defined as
\begin{equation}
\label{eq:bkw-l2-rc}
\mathcal{E}^{\mathrm{RC}}_{\Delta t_{\rm ref},r}(t)
=
\big\|\mQ_{\rm ex}(t)-\mQ_{\mathrm{RC}}^{(\Delta t_{\rm ref},r)}(f_{\rm ex}, f_{\rm ex})\big\|_2,\qquad 
\mathcal{E}^{\mathrm{FSM}}(t)
=
\big\|\mQ_{\rm ex}(t)-\mQ_{\mathrm{FSM}}(f_{\rm ex}, f_{\rm ex})\big\|_2,
\end{equation}
where $\|\cdot\|_2$ denotes the discrete $L^2$ norm on the microscopic velocity space.

\begin{figure}[!hptb]
\centering
\subfigure[Collision error]{
    \includegraphics[width=0.35\textwidth]{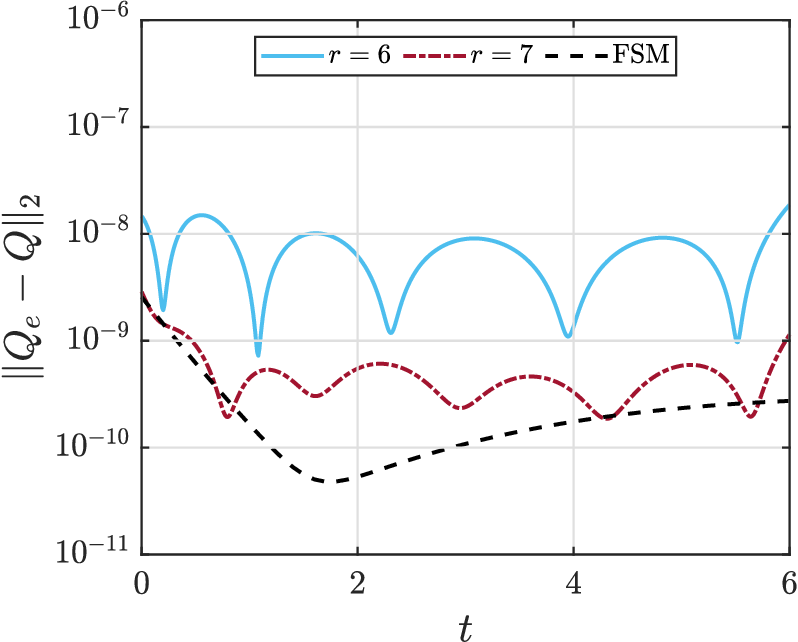}
    \label{fig:bkw_coll_err}
} \qquad 
\subfigure[Solution error]{
    \includegraphics[width=0.35\textwidth]{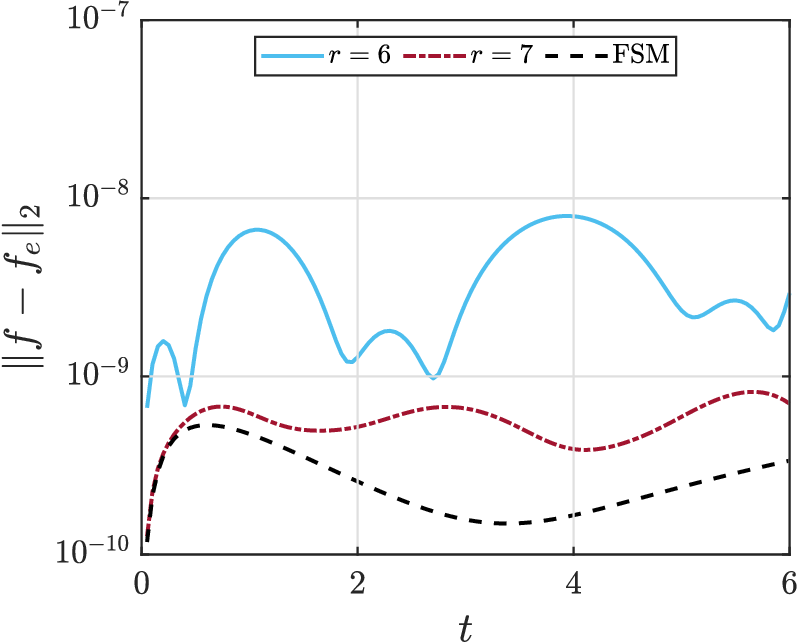}
    \label{fig:bkw_sol_err}
}
\caption{(BKW relaxation problem in Sec.~\ref{sec:BKW}) Time evolution of the $L^2$ error for the collision model and the numerical solution. (a) the collision error $\|\mQ_{\rm ex}-\mQ\|_2$ evaluated on the exact solution; (b) $L^2$ error of the numerical solution $\|f-f_{\rm ex}\|_2$, for the FSM and the RC model with $r=N_f=N_Q=6$ and $7$.}
\label{fig:bkw_error_time}
\end{figure}

The RC collision error $\mathcal{E}^{\mathrm{RC}}_{\Delta t_{\rm ref},r}$ at $t=1$ with the reference time steps $\Delta t_{\rm ref} = 0.5, 0.2, 0.1$ and $0.05$ and the ranks $r = 3, \dots, 7$ is reported in Tab. \ref{tab:Error_bkw_coll}, together with the FSM error $\mathcal{E}^{\mathrm{FSM}}$. It shows that the RC collision error decreases rapidly as the rank number increases and reaches the same order of magnitude as the FSM error at $r=7$. Moreover, the impact of the reference time step length is small. For each rank number, refining $\Delta t_{\rm ref}$ from $0.5$ to $0.05$ hardly changes the RC collision error. This is consistent with the rank selection criterion in Fig.~\ref{fig:bkw_energy}, where the decay of the residual energy is also insensitive to $\Delta t_{\rm ref}$. Therefore, the largest reference time step $\Delta t_{\rm ref} = 0.5$, which requires the smallest amount of offline data, is adopted in the rest of this section. The evolution of the collision error $\mathcal{E}^{\mathrm{RC}}_{\Delta t_{\rm ref},r}(t)$ and $\mathcal{E}^{\mathrm{FSM}}(t)$ evaluated on the exact solution $f_{\rm ex}$ is plotted in Fig. \ref{fig:bkw_coll_err}, where the rank number is set as $r = 6$ and $7$. It shows that the error with $r = 7$ is at the same order of FSM, while that with $r = 6$ is about ten times larger than that of FSM. 

\begin{table}[htbp]
\centering
{\footnotesize
\renewcommand{\arraystretch}{1.15}
\begin{tabular}{@{}c cccc@{}}
\toprule
\multirow{2}{*}{$r$} & \multicolumn{4}{c}{$\Delta t_{\rm ref}$} \\
\cmidrule(lr){2-5}
    & $0.5$ & $0.2$ & $0.1$ & $0.05$ \\
\midrule
3   & $3.63{\times}10^{-5}$  & $3.74{\times}10^{-5}$  & $3.74{\times}10^{-5}$  & $3.74{\times}10^{-5}$ \\
4   & $5.71{\times}10^{-7}$  & $5.20{\times}10^{-7}$  & $5.59{\times}10^{-7}$  & $5.68{\times}10^{-7}$ \\
5   & $1.80{\times}10^{-7}$  & $1.02{\times}10^{-7}$  & $9.61{\times}10^{-8}$  & $9.50{\times}10^{-8}$ \\
6   & $2.97{\times}10^{-9}$  & $5.50{\times}10^{-9}$  & $5.69{\times}10^{-9}$  & $5.71{\times}10^{-9}$ \\
7   & $\bf{4.47{\times}10^{-10}}$ & $2.57{\times}10^{-10}$ & $2.54{\times}10^{-10}$ & $2.56{\times}10^{-10}$ \\
\midrule
FSM & \multicolumn{4}{c}{$1.64{\times}10^{-10}$} \\
\bottomrule
\end{tabular}
}
\caption{(BKW relaxation problem in Sec.~\ref{sec:BKW}) RC collision error $\mathcal{E}^{\mathrm{RC}}_{\Delta t_{\rm ref},r}$ at $t=1$ with different reference time step length $\Delta t_{\rm ref}$ and rank numbers $r=N_f=N_Q$, together with the FSM error $\mathcal{E}^{\mathrm{FSM}}$. The offline data are generated on the entire time $[0,6]$.}
\label{tab:Error_bkw_coll}
\end{table}



\begin{figure}[!hptb]
\centering
\subfigure[$t=0.5$]{
            \includegraphics[width = 0.3\textwidth]{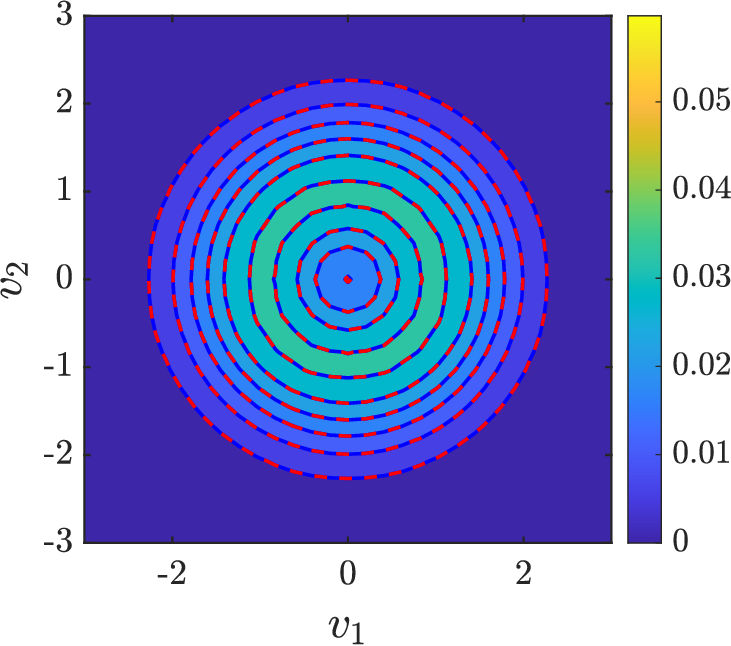}
}\hfill
\subfigure[$t=2$]{
            \includegraphics[width = 0.3\textwidth]{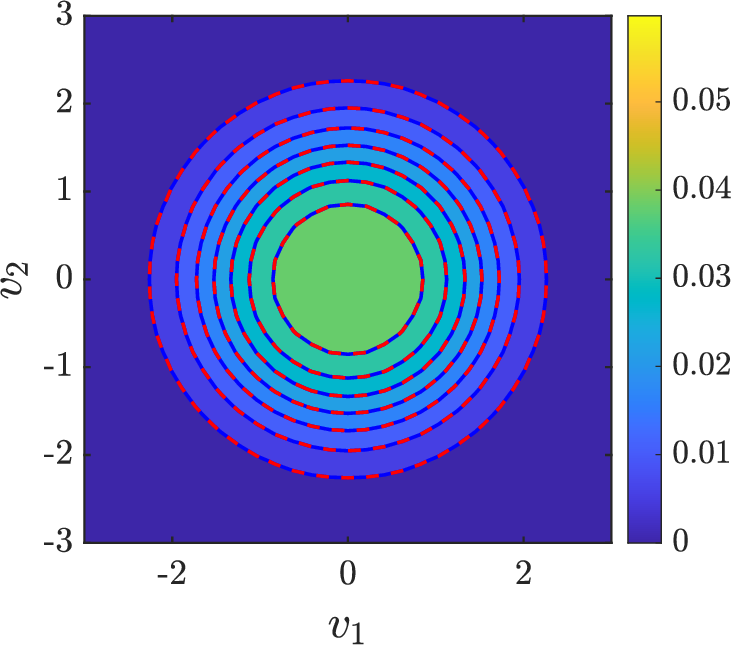}
}\hfill
\subfigure[$t=6$]{
            \includegraphics[width = 0.3\textwidth]{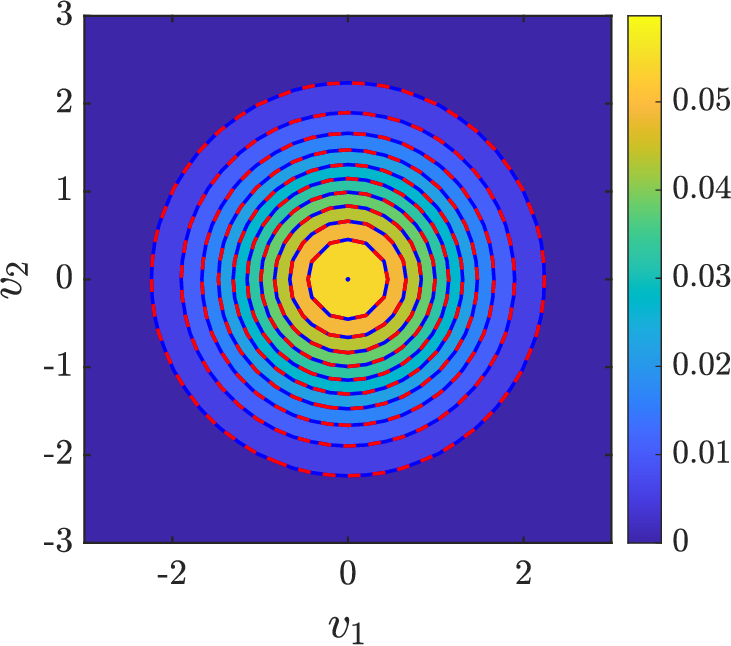}
}
\caption{(BKW relaxation problem in Sec.~\ref{sec:BKW}) Numerical solution of the distribution function on the slice $v_3 = 0$ at $t = 0.5, 2$ and $6$. Here, the solid blue lines are those of RC, the red dashed lines are those of FSM, and the filled contours are those of the exact solution. } 
\label{fig:bkw}
\end{figure}

Then, the homogeneous equation \eqref{eq:bkw} is solved numerically, where the fourth-order Runge--Kutta method with time step length $\Delta t = 0.05$ is utilized for time discretization, and the collision term is approximated by FSM and RC, respectively. The numerical solution by FSM and RC together with the exact solution at $t = 0.5, 2$ and $6$ is plotted in Fig. \ref{fig:bkw}, where the rank number for RC is set as $r = 7$, and the contours are plotted on the slice $v_3=0$ to visualize the three-dimensional distribution. It shows that the RC solution is visually indistinguishable from the FSM solution, and both of them are in close agreement with the exact solution. This shows that the reduced collision model captures the whole relaxation process from the initial non-equilibrium state to the Maxwellian equilibrium. The evolution of the $L^2$ error between the analytical solution and the numerical solution by FSM and RC is shown in Fig. \ref{fig:bkw_sol_err}, where the similar behavior of the collision error is illustrated.

\begin{table}[!hptb]
\centering
{\footnotesize
\renewcommand{\arraystretch}{1.15}
\begin{tabular*}{0.85\textwidth}{@{\extracolsep{\fill}}r r c r r l@{}}
\toprule
$N$ & method & $r$ & offline (s) & online (s) & $L^2$ error \\
\midrule
\multirow{3}{*}{8}
 & FSM & --- & --- & 18.0 & $2.7{\times}10^{-2}$ \\
\cmidrule(l){2-6}
 & \multirow{2}{*}{RC } & 1 & 1.9 & 0.013 & $5.8{\times}10^{-2}$ \\
 &                      & 2 & 2.0 & {\bf 0.013} & ${\bf 2.8{\times}10^{-2}}$ \\
\midrule
\multirow{4}{*}{12}
 & FSM & --- & --- & 101.0 & $6.4{\times}10^{-4}$ \\
\cmidrule(l){2-6}
 & \multirow{3}{*}{RC } & 1 & 10.3 & 0.028 & $5.5{\times}10^{-3}$ \\
 &                      & 2 & 11.1 & 0.037 & $4.9{\times}10^{-4}$ \\
 &                      & 3 & 12.2 & {\bf 0.046} & ${\bf 6.3{\times}10^{-4}}$ \\
\midrule
\multirow{4}{*}{16}
 & FSM & --- & --- & 153.1 & $8.0{\times}10^{-7}$ \\
\cmidrule(l){2-6}
 & \multirow{3}{*}{RC } & 2 & 16.9 & 0.047 & $2.4{\times}10^{-4}$ \\
 &                      & 3 & 18.5 & 0.049 & $8.8{\times}10^{-6}$ \\
 &                      & 4 & 20.6 & {\bf 0.060} & ${\bf 7.8{\times}10^{-7}}$ \\
\midrule
\multirow{4}{*}{24}
 & FSM & --- & --- & 1283.6 & $3.4{\times}10^{-10}$ \\
\cmidrule(l){2-6}
 & \multirow{3}{*}{RC } & 5 & 203.3 & 0.139 & $2.9{\times}10^{-8}$ \\
 &                      & 6 & 231.7 & 0.138 & $2.9{\times}10^{-9}$ \\
 &                      & 7 & 264.6 & {\bf 0.172} & ${\bf 7.0{\times}10^{-10}}$ \\
\bottomrule
\end{tabular*}
}
\caption{(BKW relaxation problem in Sec.~\ref{sec:BKW}) Efficiency comparison of RC and FSM. }
\label{tab:bkw-all-in-one}
\end{table}

Finally, the efficiency of RC is studied at the solution level. The offline data are generated by solving \eqref{eq:bkw} with the reference time step length $\Delta t_{\rm ref}=0.5$ to $t=6$, while in the online stage both the RC model and the FSM are solved with the time step length $\Delta t=0.05$ to $t=6$. 
The $L^2$ error at $t = 6$, together with the CPU time for $N = 8, 12, 16$ and $24$ and different rank numbers, is listed in Tab. \ref{tab:bkw-all-in-one}. Here, the rank number is increased until the error of RC reaches the same level as the FSM error. The computational time of the RC model is split into the offline time and the online time. The offline time consists of the generation of the offline data \eqref{eq:tr_f} and \eqref{eq:tr_Q}, the singular value decomposition \eqref{eq:svd} of the two data matrices, and the construction of the kernel tensor \eqref{eq:kerner-tensor}, while the online time is the cost of solving \eqref{eq:bkw} with the reduced collision model. 

Tab. \ref{tab:bkw-all-in-one} shows that a larger $N$ requires a larger rank, since the FSM solution itself becomes more accurate when the microscopic velocity space is refined. The offline time is smaller than the CPU time of the FSM for all the cases reported here, since the offline data are generated with the much larger time step $\Delta t_{\rm ref}$. As a result, the CPU time for the RC model is five to nine times shorter than that of FSM, even with the offline cost included. In particular, the CPU time for one evaluation of the reduced collision model is almost negligible, and the online time of RC is more than three orders of magnitude smaller than that of FSM.

\subsection{Steady-state Boltzmann equation}
\label{subsec:steady}
Next, the steady-state Boltzmann equation is considered, where the spatially 1D channel problems are studied in Sec. \ref{sec:fourier} and \ref{sec:couette}, and the spatially 2D problem is tested for the lid-driven cavity problem in Sec. \ref{subsubsec:cavity}. For the one-dimensional channel problems, the physical domain is $\Omega=[0,1]$, it reduces to
\begin{equation}
v_1\,\partial_x f(x,\bv)=\frac{1}{\Kn}\,\mQ(f,f)(x,\bv).
\label{eq:steady_1d}
\end{equation}
The source iteration (SI) is utilized to solve this steady-state problem \cite{valougeorgis2003acceleration, lihnaropoulos2007formulation}. The spatial domain is discretized by a uniform mesh $\{x_j\}_{j=1}^{N_x}$ with the mesh size $\Delta x$, while the microscopic velocity space follows the DVM discretization described in Sec.~\ref{sec:dvm}. Let $f_{j,s}^{(n)}$ denote the numerical solution of $f(x_j,\bv_s)$ at the $n$-th SI step. The numerical scheme has the form below 
\begin{equation}
v_{1,s}^{+}\frac{f_{j,s}^{(n+1)}-f_{j-1,s}^{(n+1)}}{\Delta x}
+
v_{1,s}^{-}\frac{f_{j+1,s}^{(n+1)}-f_{j,s}^{(n+1)}}{\Delta x}
+\nu f_{j,s}^{(n+1)}
=
\frac{1}{\Kn}\mQ_{j,s}^{(n)}+\nu f_{j,s}^{(n)}.
\label{eq:si}
\end{equation}
Here, $v_{1,s}^{\pm}=(v_{1,s}\pm|v_{1,s}|)/2$, and $\nu$ is often chosen as the collision frequency, which depends on the density and temperature of the gas. In each SI sweep, the dominant cost lies in the evaluation of the quadratic collision term $\mQ(f^{(n)},f^{(n)})$.

For the boundary condition, the Maxwell wall boundary condition is imposed \cite{shen2006rarefied}. Precisely, at a wall point $\bx_b$ with wall velocity $\bu^W$, wall temperature $T^W$, and inward unit normal $\bn$, the outgoing distribution is prescribed as
\begin{equation}
f(\bx_b,\bv)
=
\frac{\rho^W}{(2\pi T^W)^{3/2}}
\exp\left(-\frac{|\bv-\bu^W|^2}{2T^W}\right),
\qquad \bc^W\cdot\bn>0.
\label{eq:boun_con}
\end{equation}
Here, $\bc^W=\bv-\bu^W$ is the peculiar velocity relative to the wall. The wall density $\rho^W$ is determined from the impermeability condition
\begin{equation}
    \label{eq:flux}
    \int_{\bbR^3}(\boldsymbol{c}^W \cdot \boldsymbol{n}) f(\bx_b,\bv) \dd\bv = 0.
\end{equation}
Substituting \eqref{eq:boun_con} into \eqref{eq:flux} yields
\begin{equation}
\label{eq:bound_rho}
    \rho^W = -\dfrac{\int_{\boldsymbol{c}^W \cdot \boldsymbol{n}<0}(\boldsymbol{c}^W \cdot \boldsymbol{n}) f_{\text{in}}(\bx_b,\bv) \dd\bv}{\int_{\boldsymbol{c}^W \cdot \boldsymbol{n}>0} \frac{
(\boldsymbol{c}^W \cdot \boldsymbol{n})    
    }{(2\pi T^W)^{3/2}}\exp\left(-\frac{|\bv-\boldsymbol{u}^W|^2}{2T^W}\right) \dd\bv},
\end{equation} 
where $f_{\mathrm{in}}$ denotes the incoming distribution obtained from the interior solution.

For the channel problems, convergence is monitored through the discrete $L^2$ norm of the residual. The residual at the $n$-th SI step is defined as
\begin{equation}
    r_{j,s}^{(n)} = v_{1,s}^{+} \frac{f_{j,s}^{(n)} - f_{j-1,s}^{(n)}}{\Delta x} + v_{1,s}^{-} \frac{f_{j+1,s}^{(n)} - f_{j,s}^{(n)}}{\Delta x} - \frac{1}{\Kn} \mQ_{j,s}^{(n)},
\end{equation}
with residual norm
\begin{equation}
\|\mathrm{Res}^{(n)}\|_2
=
\left(
\sum_{j,s}
\left(\frac{L}{N}\right)^3
\Delta x\,|r_{j,s}^{(n)}|^2
\right)^{1/2}.
\label{eq:residual_norm_1d}
\end{equation}
The iteration is terminated when $\|\mathrm{Res}^{(n)}\|_2 < 10^{-5}$. For the 2D cavity problem, the same convergence criterion is applied to the residual associated with the 2D SI.

\subsubsection{Fourier flow}
\label{sec:fourier}
The steady-state Fourier flow problems with different wall temperatures are studied in this section. Similar problems are also studied in \cite{jaiswal2019discontinuous}. Maxwell's wall boundary conditions are imposed at $x=0$ and $x=1$, with wall temperatures $T_l^W$ and $T_r^W$. Three pairs of wall temperatures are considered, which are listed in Tab.~\ref{tab:fourier_bc}. The Knudsen number is taken as $\Kn=1, 0.1$, and $0.01$. In all simulations, the number of the Fourier modes is set as $N = 16$. 
\begin{table}[!htbp]
\centering
{\footnotesize
\renewcommand{\arraystretch}{1.5}
\begin{tabular}{@{}ccc@{}}
\toprule
case & $T_l^W$ & $T_r^W$ \\
\midrule
(\text{i})   & $(T_0-10)/T_0$  & $(T_0+10)/T_0$ \\
(\text{ii})  & $(T_0-100)/T_0$ & $(T_0+100)/T_0$ \\
(\text{iii}) & $1$             & $3$ \\
\bottomrule
\end{tabular}}
\caption{(1D Fourier flow problem in Sec.~\ref{sec:fourier}) Three pairs of wall temperatures with $T_0=273.15$.}
\label{tab:fourier_bc}
\end{table}

The reduced approximation is first tested at the collision operator level, which can instruct how to choose the rank number $r$ and the resolution in spatial space to obtain the reference data when building the reduced collision model. Take case (ii) in Tab. \ref{tab:fourier_bc} as an example. The reduced collision operator $\mQ_{\mathrm{RC},r}^{(N_x)}$ is constructed from reference data generated on a spatial mesh with $N_x$ cells at rank $N_f=N_Q=r$. Define the collision error of RC as 
\begin{equation}
  \mathcal{E}_{N_x,r}
  \;=\;
  \big\|
  \mQ_{\mathrm{FSM}}(f^{\rm st},f^{\rm st})-\mQ_{\mathrm{RC},r}^{(N_x)}(f^{\rm st},f^{\rm st})
  \big\|_2,
  \label{eq:fourier-gap}
\end{equation}
where $\|\cdot\|_2$ denotes the discrete $L^2$ norm over the phase space with the weight $\Delta x\,(L/N)^3$, and $f^{\rm st}$ is the steady solution obtained by the FSM-based source iteration \eqref{eq:si}. The error $ \mathcal{E}_{N_x,r}$ with $N_x = 10, 20, 40, 80$ and the rank $r=10,20,30,40$ is listed in Tab. \ref{tab:Error_Fourier}. The results show that the collision error of RC decreases rapidly as the rank increases. For the larger ranks, the error also decreases as the number of the spatial cells $N_x$ increases, but at a considerably slower rate. These numerical results indicate that reference data generated on a coarse spatial mesh already yield a reduced collision model accurate enough for this inhomogeneous problem, which also provides a numerical basis for the following tests of this Fourier flow problem.

\begin{table}[htbp]
\centering
{\footnotesize
\renewcommand{\arraystretch}{1.15}
\begin{tabular}{@{}c cccc@{}}
\toprule
\multirow{2}{*}{$r$} & \multicolumn{4}{c}{$N_x$} \\
\cmidrule(lr){2-5}
    & $10$ & $20$ & $40$ & $80$ \\
\midrule
10  & $1.44{\times}10^{-4}$ & $1.52{\times}10^{-4}$ & $1.44{\times}10^{-4}$ & $1.37{\times}10^{-4}$ \\
20  & $9.98{\times}10^{-6}$ & $8.24{\times}10^{-6}$ & $6.31{\times}10^{-6}$ & $3.63{\times}10^{-6}$ \\
30  & $5.04{\times}10^{-6}$ & $3.61{\times}10^{-6}$ & $2.51{\times}10^{-6}$ & $9.50{\times}10^{-7}$ \\
40  & $\bf{1.35{\times}10^{-6}}$ & $9.11{\times}10^{-7}$ & $8.28{\times}10^{-7}$ & $3.25{\times}10^{-7}$ \\
\bottomrule
\end{tabular}}
\caption{(1D Fourier flow problem in Sec.~\ref{sec:fourier}) The collision error of RC $\mathcal{E}_{N_x,r}$ for case (ii) in Tab. \ref{tab:fourier_bc} at $\Kn=0.1$ with different  $N_x$ and rank $r=N_f=N_Q$.}
\label{tab:Error_Fourier}
\end{table}

\begin{figure}[!hptb]
\centering
\subfigure[$\rho$, case (\text{i})]{
            \includegraphics[width = 0.3\textwidth]{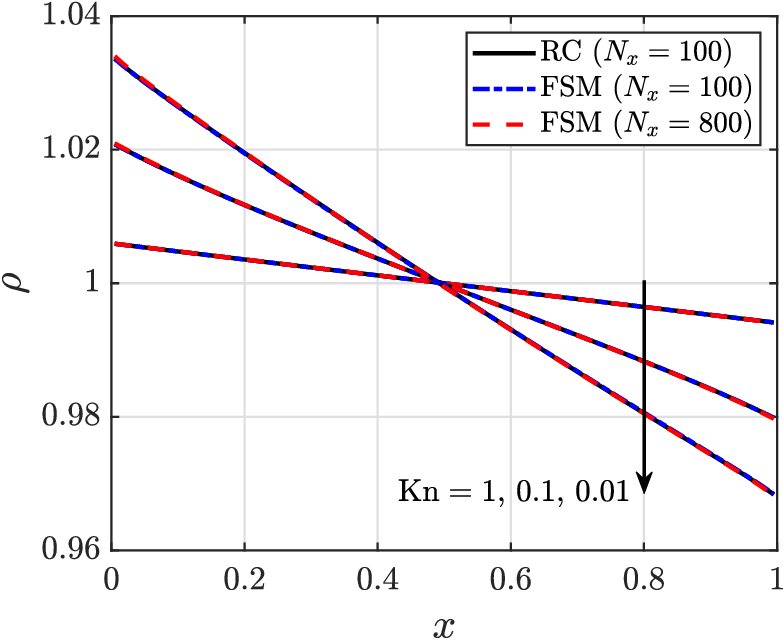}
}
\subfigure[$\rho$, case (\text{ii})]{
            \includegraphics[width = 0.3\textwidth]{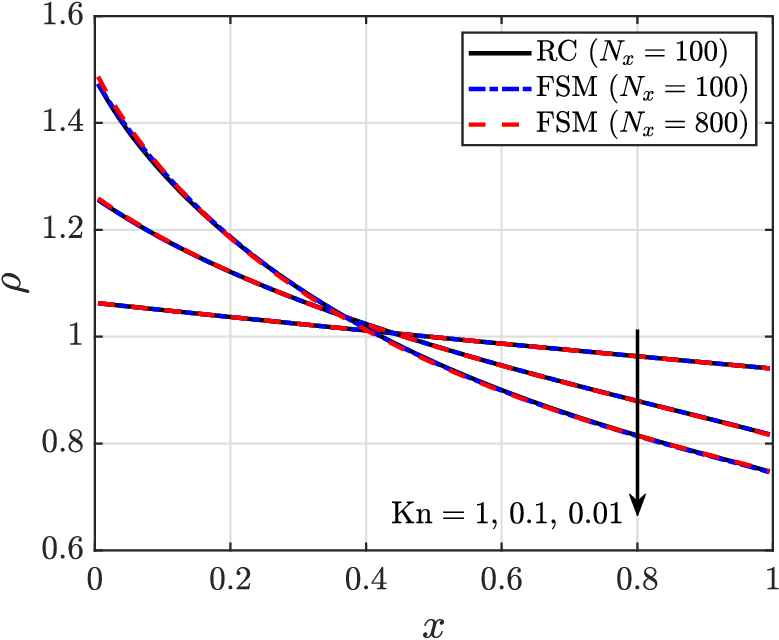}
}
\subfigure[$\rho$, case (\text{iii})]{
            \includegraphics[width = 0.3\textwidth]{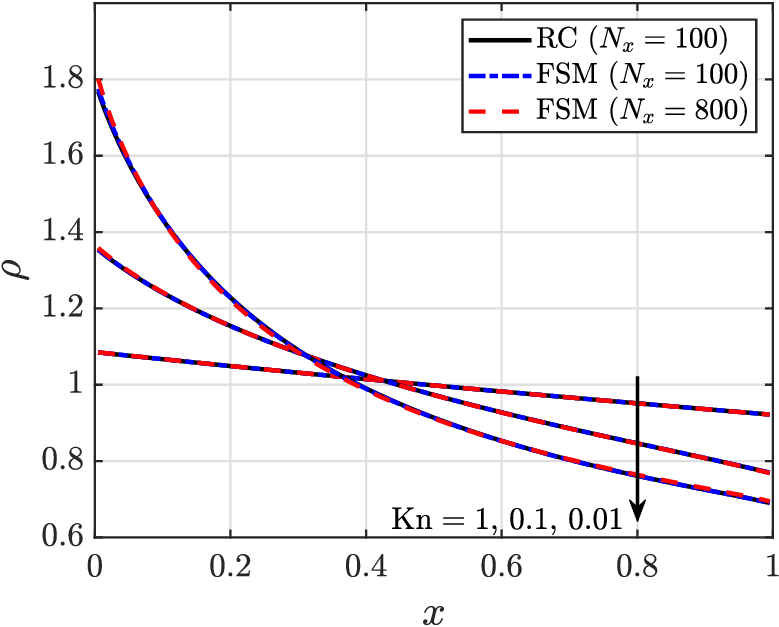}
}
\subfigure[$T$, case (\text{i})]{
            \includegraphics[width = 0.3\textwidth]{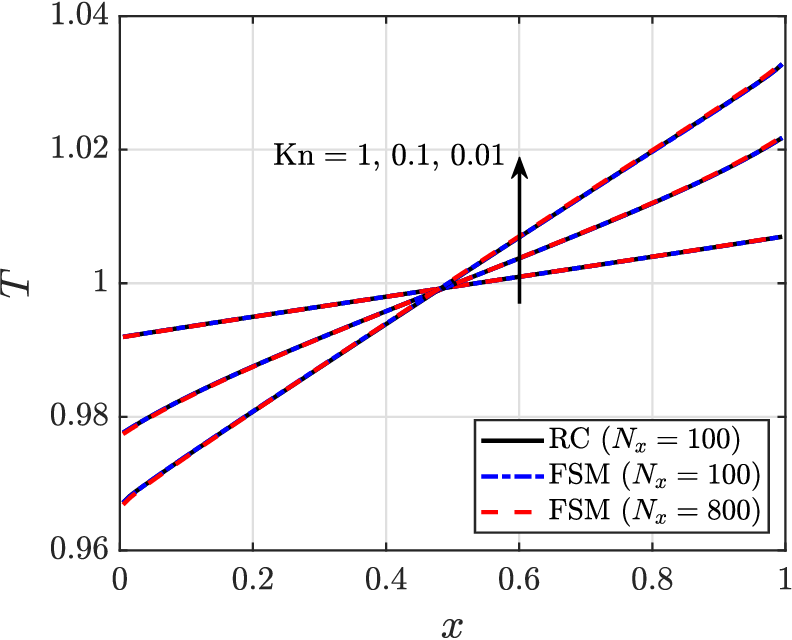}
}
\subfigure[$T$, case (\text{ii})]{
            \includegraphics[width = 0.3\textwidth]{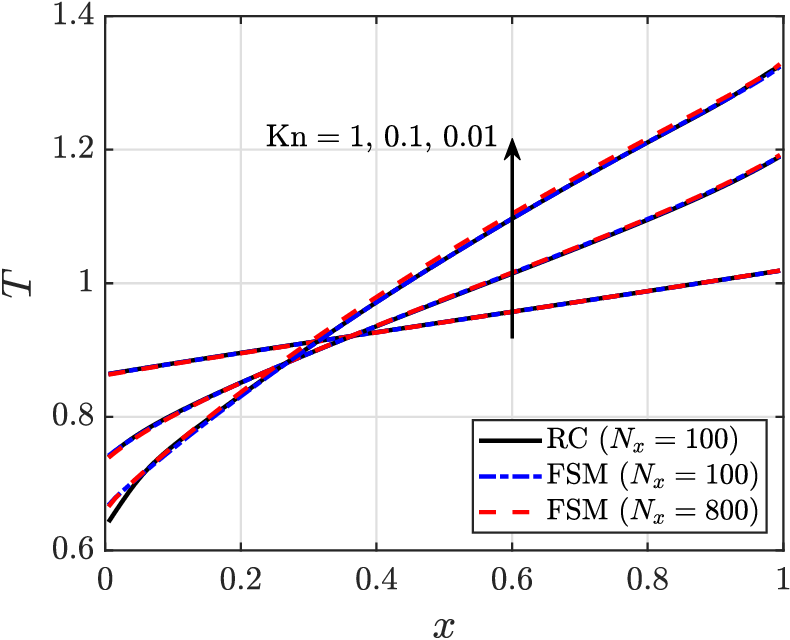}
}
\subfigure[$T$, case (\text{iii})]{
            \includegraphics[width = 0.3\textwidth]{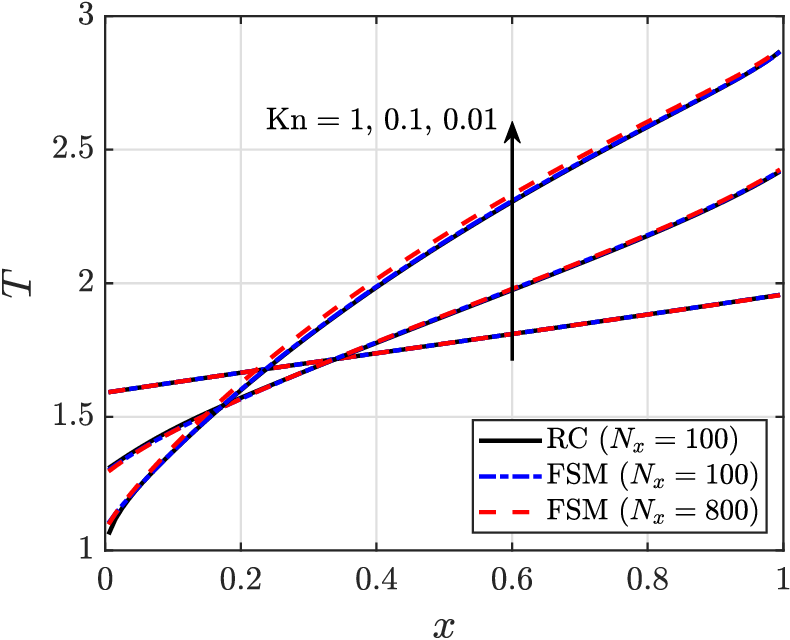}
}
\caption{
(1D Fourier flow problem in Sec.~\ref{sec:fourier}) Numerical solution of the density $\rho$ and temperature $T$ for three cases of the wall-temperatures in Tab. \ref{tab:fourier_bc} with $\Kn=1,0.1,0.01$.}
\label{fig:fourier}
\end{figure}

Next, the performance of the RC approximation is studied when solving the Fourier flow problem. The number of the spatial cells is set as $N_x = 100$ for both RC and FSM. For the RC model, it is constructed using offline data generated with rank $r = N_f = N_Q = 40$ and $N_x = 10$. The numerical solution of the density $\rho$ and temperature $T$ for the three cases of wall temperatures listed in Tab. \ref{tab:fourier_bc} at $\Kn = 1, 0.1$ and $0.01$ is shown in Fig. \ref{fig:fourier}, where the reference solution is obtained by FSM with $N_x = 800$. It shows that for all three cases, the numerical solution matches well with the reference solution obtained by FSM either with $N_x = 100$ or $N_x = 800$, from the transition regime at $\Kn=1$ to the slip and near-continuum regimes at $\Kn=0.1$ and $0.01$.

Then, the computational efficiency of the RC approximation is examined. For a macroscopic quantity, the $L^2$ error is measured by the discrete $L^2$ norm over the spatial mesh, where the numerical solution obtained by RC or FSM with $N_x=100$ is compared with the reference solution obtained by FSM with $N_x=800$. The case (ii) in Tab.~\ref{tab:fourier_bc} is studied, where the $L^2$ errors of the density and the temperature for both RC and FSM are listed in Tab.~\ref{tab:fourier-all-in-one}, together with the iteration steps and the total CPU time. Here, the CPU time for RC includes the generation of the offline data, SVD to obtain the basis function, the construction of the RC model, and the online source iteration. The detailed CPU time for different Knudsen number is listed in Tab. \ref{tab:fourier-time}. It shows that as $\Kn$ decreases, the number of SI steps increases rapidly and the total CPU cost grows accordingly. However, the iteration steps for FSM and RC are almost the same for all Knudsen numbers, and the $L^2$ errors of the density and the temperature are also of the same order. This means that the reduced collision model preserves the convergence behavior of the source iteration and does not affect the accuracy of the numerical solution. Tab. \ref{tab:fourier-all-in-one} also shows that the CPU time for RC is much shorter than FSM for $\Kn = 0.1$ and $0.01$, indicating the high efficiency of RC. 



\begin{table}[!hptb]
\centering
\begin{threeparttable}
{\footnotesize
\renewcommand{\arraystretch}{1.1}
\begin{tabular*}{0.85\textwidth}{@{\extracolsep{\fill}}r r r r c c@{}}
\toprule
\multirow{2}{*}{$\Kn$} & \multirow{2}{*}{method} & \multirow{2}{*}{steps}
& \multirow{2}{*}{CPU (s)} & \multicolumn{2}{c}{$L^2$ error} \\
\cmidrule(lr){5-6}
& & & & $\rho$ & $T$ \\
\midrule
\multirow{2}{*}{1}
& FSM & 16   & 1225  & $6.81{\times}10^{-6}$ & $1.68{\times}10^{-5}$ \\
& RC  & 16   & 1381  & $6.82{\times}10^{-6}$ & $1.44{\times}10^{-5}$ \\
\addlinespace
\multirow{2}{*}{0.1}
& FSM & 84   & 6811  & $4.65{\times}10^{-5}$ & $6.92{\times}10^{-5}$ \\
& RC  & 85   & 2007  & $4.56{\times}10^{-5}$ & $8.54{\times}10^{-5}$ \\
\addlinespace
\multirow{2}{*}{0.01}
& FSM & 1096 & 87584 & $1.41{\times}10^{-4}$ & $1.49{\times}10^{-4}$ \\
& RC  & 1099 & 12103 & $9.84{\times}10^{-5}$ & $1.40{\times}10^{-4}$ \\
\bottomrule
\end{tabular*}
}
\end{threeparttable}
\caption{(1D Fourier flow problem in Sec.~\ref{sec:fourier}) The comparison of efficiency between RC and FSM for case (\text{ii}) in Tab. \ref{tab:fourier_bc}. Here, the CPU time for RC includes both the offline construction of the RC model and the online source iteration. }
\label{tab:fourier-all-in-one}
\end{table}

\begin{table}[!hptb]
\centering
{\footnotesize
\renewcommand{\arraystretch}{1.1}
\begin{tabular*}{0.75\textwidth}{@{\extracolsep{\fill}}r r c r r@{}}
\toprule
\multirow{2}{*}{$\Kn$} & \multicolumn{2}{c}{offline (s)} & \multirow{2}{*}{online SI (s)} & \multirow{2}{*}{total (s)} \\
\cmidrule(lr){2-3}
 & data generation & SVD and kernel & & \\
\midrule
1    & 112  & 1263 & 6   & 1381 \\
0.1  & 672  & 1303 & 32  & 2007 \\
0.01 & 8792 & 2880 & 431 & 12103 \\
\bottomrule
\end{tabular*}
}
\caption{(1D Fourier flow problem in Sec.~\ref{sec:fourier}) Detailed CPU time of the RC model in Tab.~\ref{tab:fourier-all-in-one} for case (\text{ii}) in Tab. \ref{tab:fourier_bc}.}
\label{tab:fourier-time}
\end{table}

\subsubsection{Couette flow}
\label{sec:couette}
In this section, the steady-state Couette flow problem between two parallel isothermal walls moving with opposite tangential velocities is considered. Maxwell's wall boundary conditions are imposed at $x=0$ and $x=1$, with wall velocities
\begin{equation}
    \bu_l^W=(0,-U^W,0),\qquad \bu_r^W=(0, U^W,0),
\end{equation}
while the two walls are kept at the same temperature. Two cases $U^W = 0.5$ and $U^W = 1$ with three Knudsen numbers $\Kn = 1, 0.1$ and $0.01$ are tested. A similar example is also tested in \cite{shakhov1969couette, jaiswal2019discontinuous}.

The numerical setup is the same as that utilized for the Fourier flow problem in Sec.~\ref{sec:fourier}. The number of Fourier modes is set as $N=16$, and the number of the spatial cells is $N_x = 100$ for both RC and FSM. The reference solution is obtained by FSM with $N_x=800$, and the rank number $N_f = N_Q = 40$ is utilized in RC, with reduced bases constructed offline from offline data generated on the coarse grid $N_x=10$.

\begin{figure}[!hptb]
\centering
\subfigure[$\rho$, $U^W=0.5$]{
            \includegraphics[width = 0.225\textwidth]{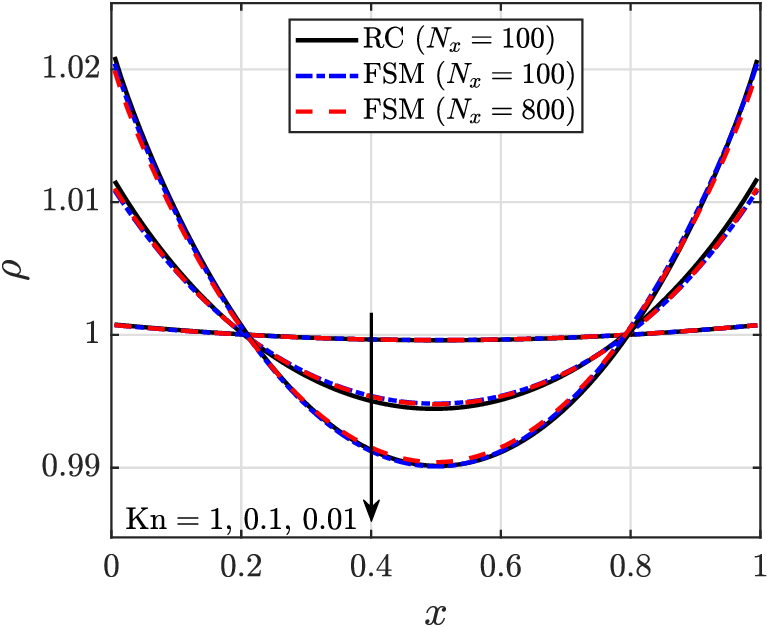}
}\subfigure[$u_y$, $U^W=0.5$]{
            \includegraphics[width = 0.23\textwidth]{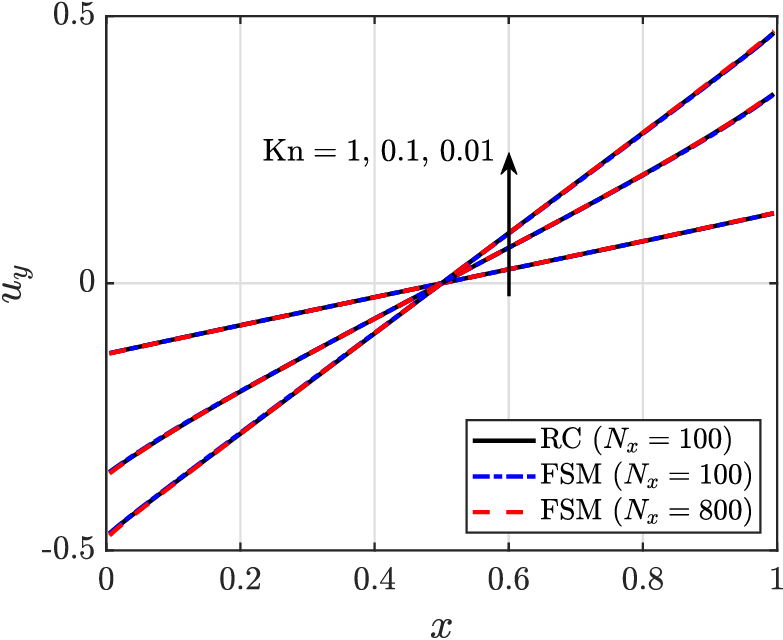}
}
\subfigure[$\rho$, $U^W=1$]{
            \includegraphics[width = 0.23\textwidth]{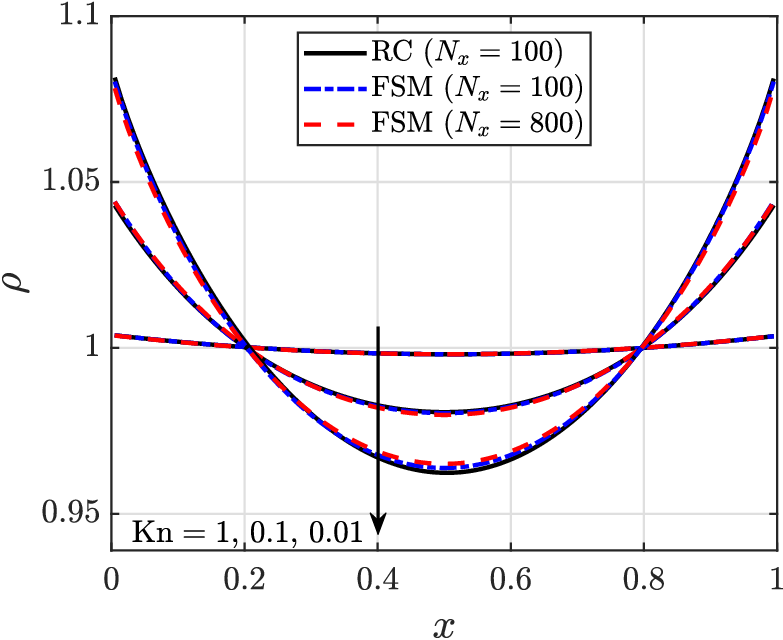}
}
\subfigure[$u_y$, $U^W=1$]{
            \includegraphics[width = 0.23\textwidth]{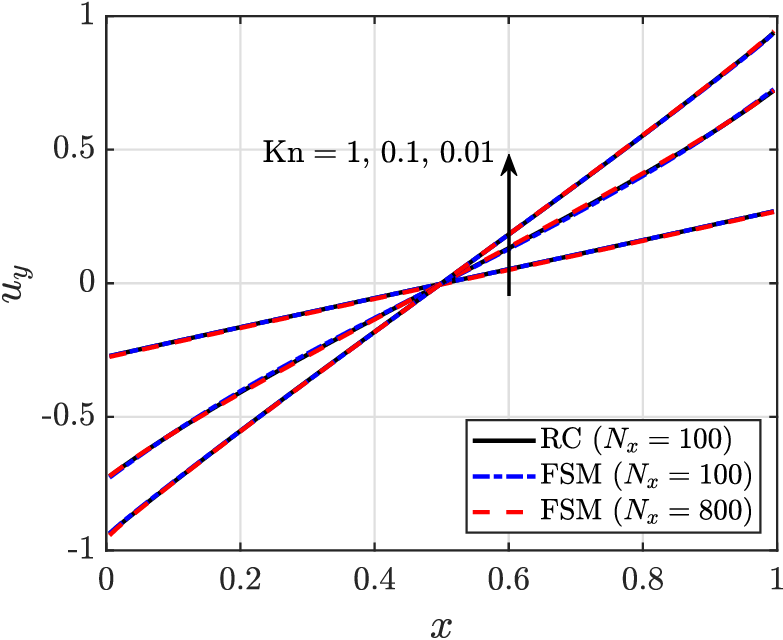}
}
\caption{
(1D Couette flow problem in Sec.~\ref{sec:couette}) Steady-state numerical solution of the density $\rho$ and the tangential velocity $u_y$ for the two wall speeds $U^W=0.5$ and $U^W=1$ at $\Kn=1,0.1$ and $0.01$.}
\label{fig:couette}
\end{figure}

The steady-state solution of the density $\rho$ and the tangential velocity $u_y$ with the two wall speeds for $\Kn = 1, 0.1$ and $0.01$ is plotted in Fig. \ref{fig:couette}. For all the cases, the numerical solution of RC is nearly indistinguishable from that of FSM on the same mesh, and remains in close agreement with the reference solution obtained with FSM with $N_x = 800$. Similar to the Fourier flow problem, the agreement between RC and FSM confirms the accuracy of RC in capturing rarefaction phenomena such as wall slip in steady one-dimensional problems.

\subsubsection{2D lid-driven cavity flow}
\label{subsubsec:cavity}
Next, the classical two-dimensional lid-driven cavity flow problem in the unit square $\Omega=(0,1)\times(0,1)$ is considered, for which the steady-state Boltzmann equation  becomes
\begin{equation}
v_1\partial_{x}f+v_2\partial_{y}f=\frac{1}{\Kn}\,\mQ(f,f).
\end{equation}
Here, the Maxwell's wall boundary conditions are imposed at all walls with the wall temperature $T^W = 1$ and wall velocity $\bu^W = \bm{0}$, except for the top wall whose wall velocity is 
\begin{equation}
    \bu^W=(0.5,0,0).
\end{equation}
The Knudsen number is taken as $\Kn=1,0.1$, and $0.01$. Similar problems are also studied in \cite{naris2005driven, jaiswal2019discontinuous}.

The same source iteration as the 1D problems is utilized here, with the first-order upwind numerical flux applied in both spatial directions. The number of the Fourier modes is set as $N = 16$, and the spatial domain is discretized on a uniform mesh with $N_x = N_y = 100$ for both RC and FSM. The rank number is increased to $N_f = N_Q = 80$ to account for the higher complexity of the two-dimensional dynamics, with the reduced bases constructed from offline data generated on the coarse mesh $N_x = N_y = 10$.

\begin{figure}[!hptb]
\centering
\subfigure[$\rho$]{
            \includegraphics[width = 0.23\textwidth]{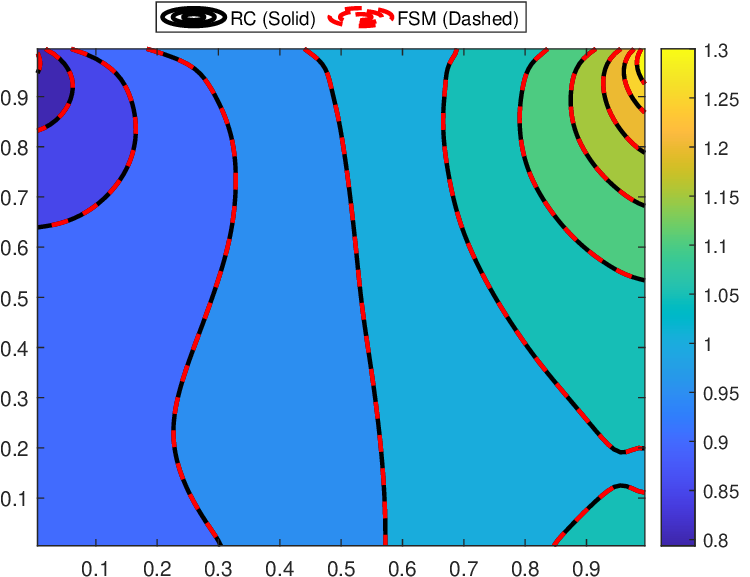}
}
\subfigure[$u_x$]{
            \includegraphics[width = 0.23\textwidth]{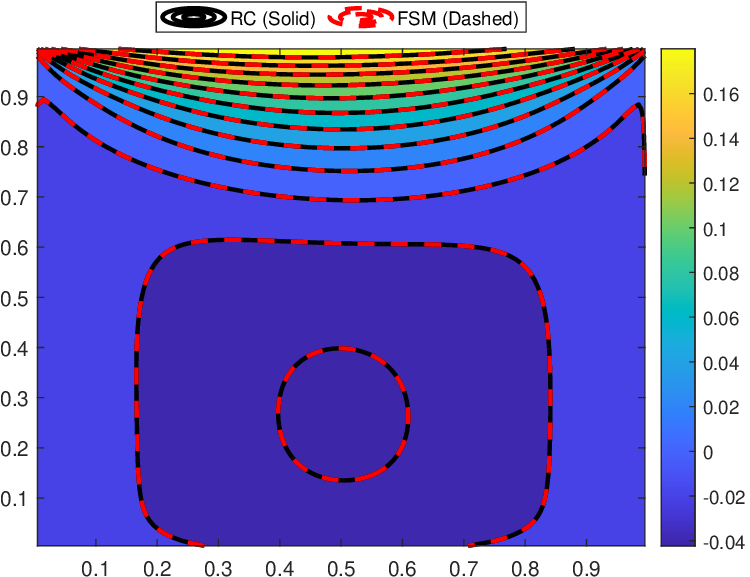}
}
\subfigure[$u_y$]{
            \includegraphics[width = 0.23\textwidth]{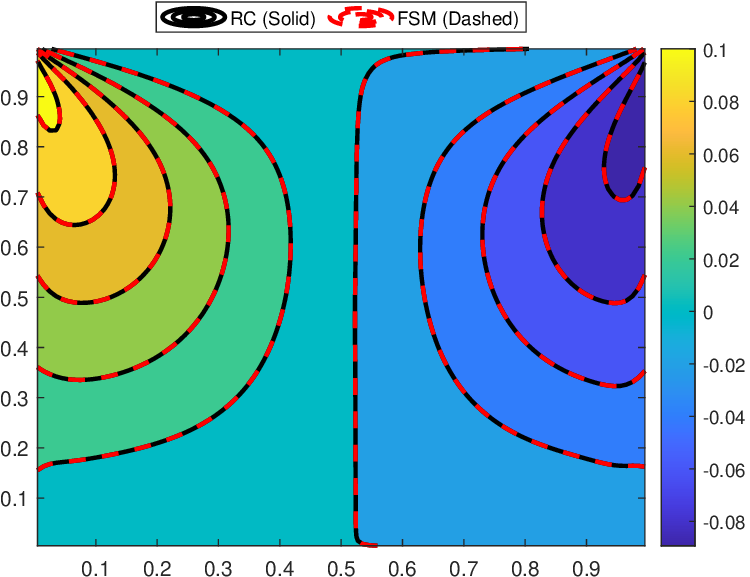}
}
\subfigure[$T$]{
            \includegraphics[width = 0.23\textwidth]{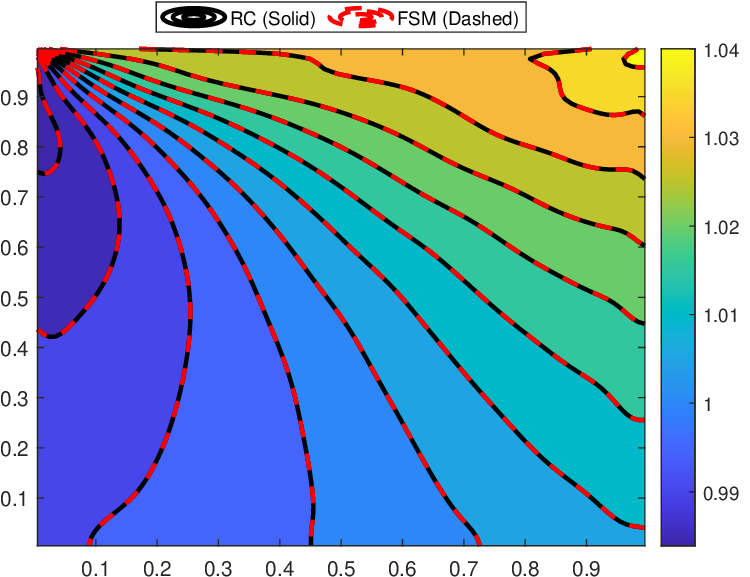}
}
\caption{(2D lid-driven cavity flow problem in Sec.~\ref{subsubsec:cavity}) Steady-state numerical solution of the density $\rho$, the macroscopic velocity $u_x$, $u_y$ and the temperature $T$ for $\Kn = 1$. Here, the solid black lines are those of RC and the red dashed lines are those of FSM. }
\label{fig:Cavity_kn1}
\end{figure}

\begin{figure}[!hptb]
\centering
\subfigure[$\rho$]{
            \includegraphics[width = 0.23\textwidth]{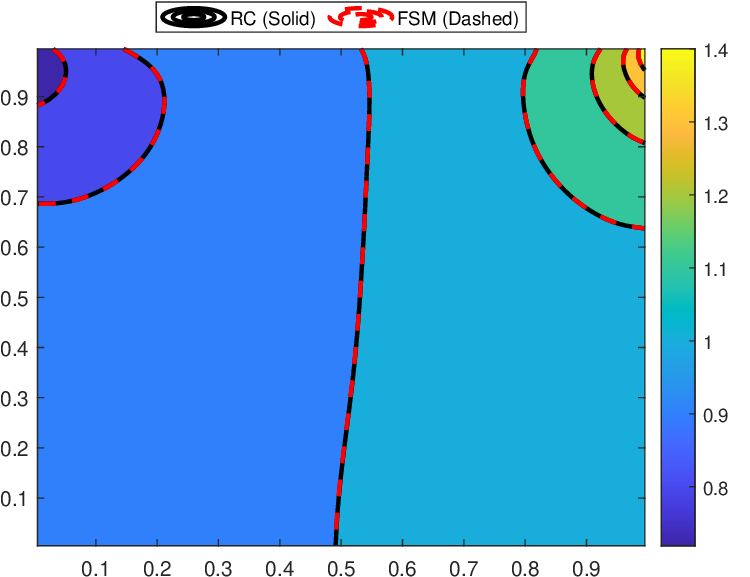}
}
\subfigure[$u_x$]{
            \includegraphics[width = 0.23\textwidth]{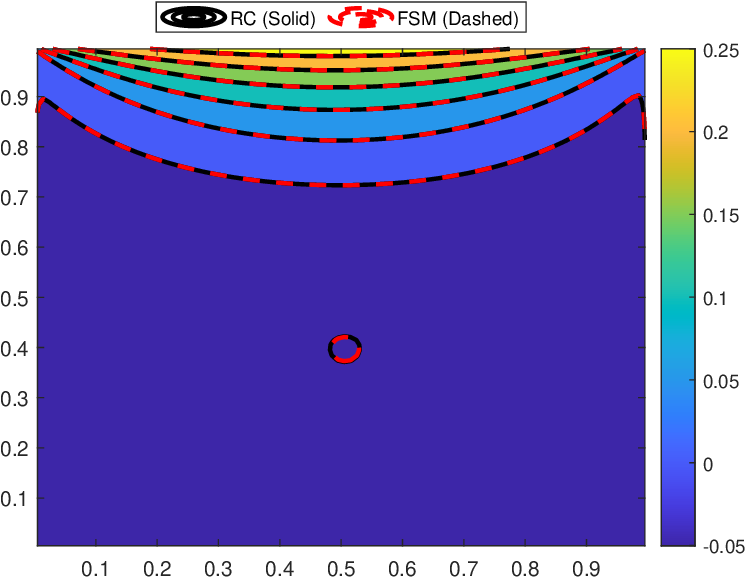}
}
\subfigure[$u_y$]{
            \includegraphics[width = 0.23\textwidth]{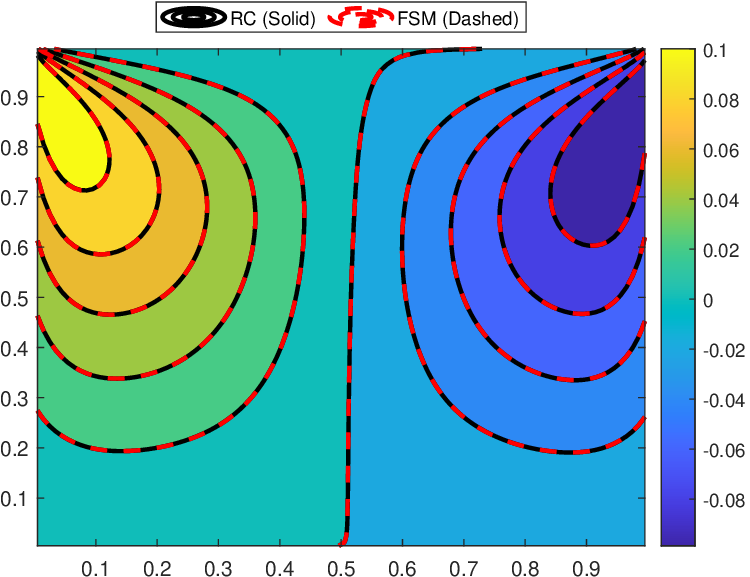}
}
\subfigure[$T$]{
            \includegraphics[width = 0.23\textwidth]{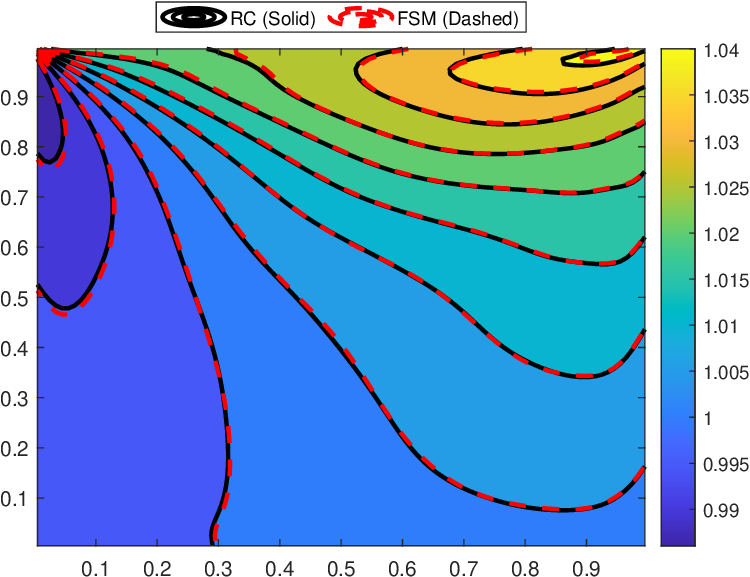}
}
\caption{(2D lid-driven cavity flow problem in Sec.~\ref{subsubsec:cavity}) Steady-state numerical solution of the density $\rho$, the macroscopic velocity $u_x$, $u_y$ and the temperature $T$ for $\Kn = 0.1$. Here, the solid black lines are those of RC and the red dashed lines are those of FSM.
}
\label{fig:Cavity_kndot1}
\end{figure}

\begin{figure}[!hptb]
\centering
\subfigure[$\rho$]{
            \includegraphics[width = 0.23\textwidth]{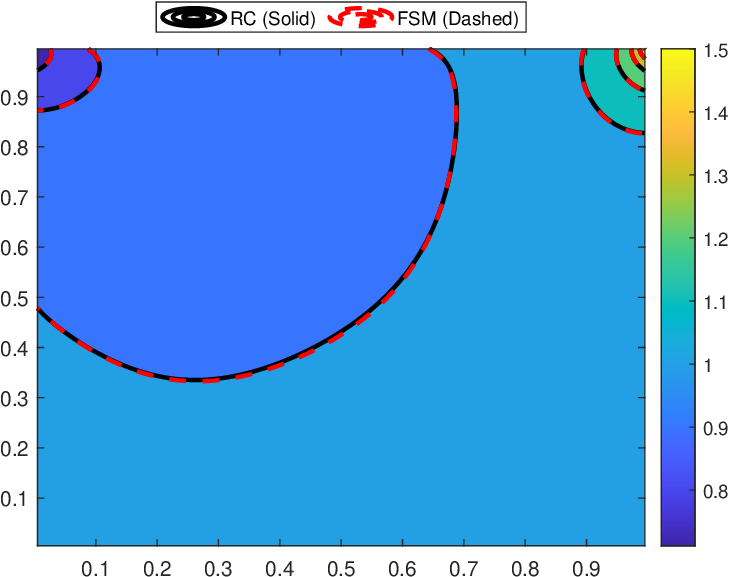}
}
\subfigure[$u_x$]{
            \includegraphics[width = 0.23\textwidth]{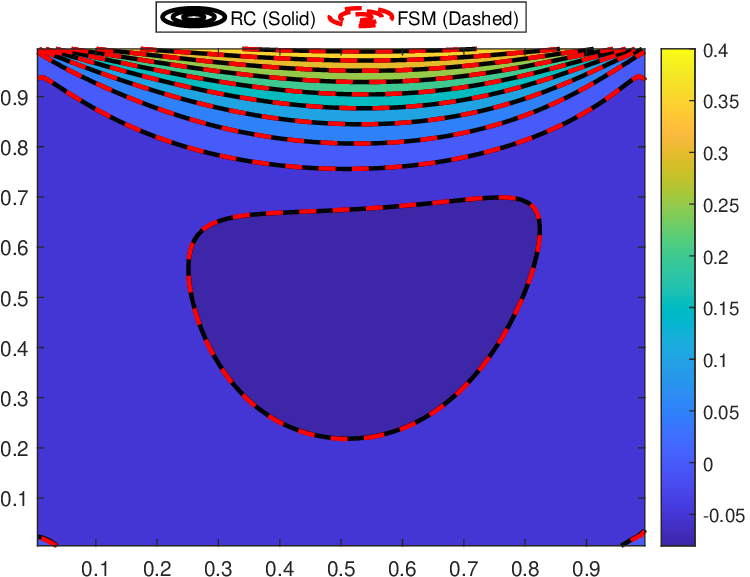}
}
\subfigure[$u_y$]{
            \includegraphics[width = 0.23\textwidth]{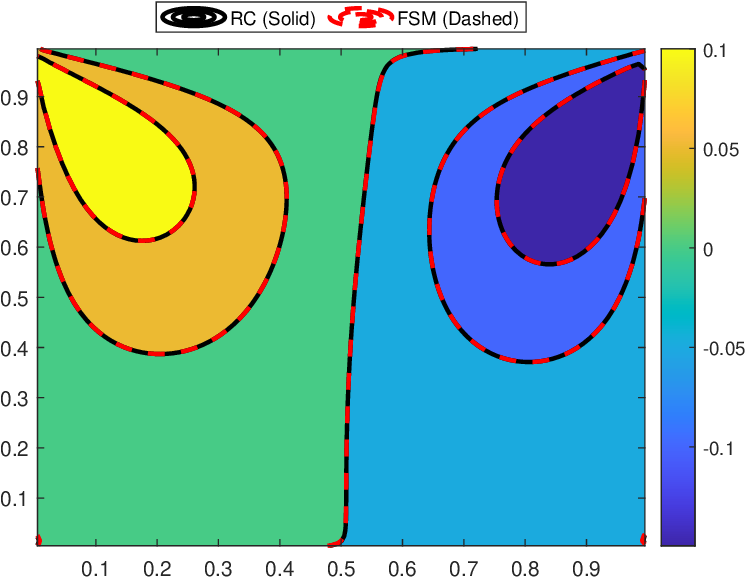}
}
\subfigure[$T$]{
            \includegraphics[width = 0.23\textwidth]{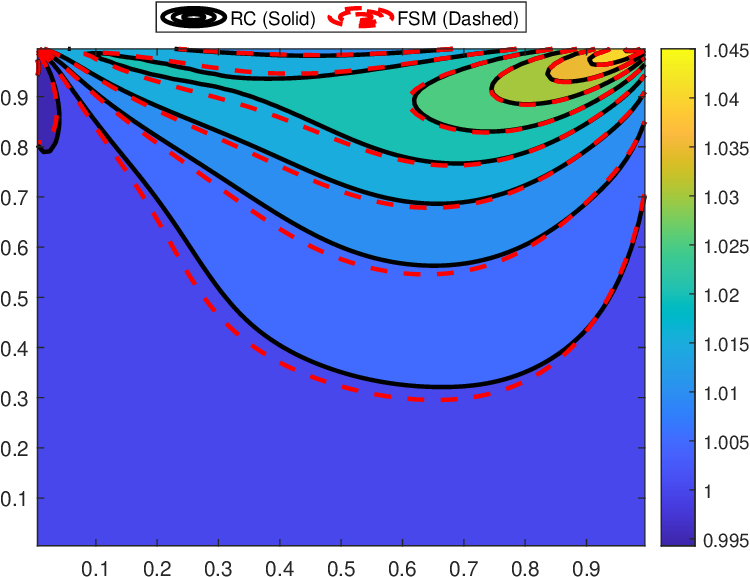}
}
\caption{(2D lid-driven cavity flow problem in Sec.~\ref{subsubsec:cavity}) 
Steady-state numerical solution of the density $\rho$, the macroscopic velocity $u_x$, $u_y$ and the temperature $T$ for $\Kn = 0.01$. Here, the solid black lines are those of RC and the red dashed lines are those of FSM.}
\label{fig:Cavity_kndot01}
\end{figure}

The steady state numerical solution of the density $\rho$, the macroscopic velocity $u_x$, $u_y$ and the temperature $T$ for $\Kn = 1, 0.1$ and $0.01$ is plotted in Fig. \ref{fig:Cavity_kn1} to Fig. \ref{fig:Cavity_kndot01}, respectively. For all $\Kn$, they show that the numerical solution by RC is almost the same as the reference solution obtained by FSM. Moreover, as $\Kn$ decreases, the gas near the lid is dragged closer to the lid speed, so that the velocity slip at the moving wall becomes weaker. At the same time, the transverse velocity $u_y$ shows that the circulation induced by the lid is confined to a region closer to the lid, while the lower part of the cavity becomes nearly quiescent. Both trends are well reproduced by the RC model.

Then, the computational efficiency of RC is examined. A reference solution with a finer mesh is prohibitively expensive for this 2D problem.
Therefore, the $L^2$ error is measured by the discrete $L^2$ norm over the spatial mesh of the difference between the RC and the FSM solutions obtained on the same mesh with $N_x = N_y = 100$. The errors of the density, the macroscopic velocity and the temperature are listed in Tab.~\ref{tab:cavity_time}, where the number of the SI steps and the CPU time are also given. It shows that for all $\Kn$, these $L^2$ errors are all of the same order, namely $10^{-5}$ for the density and the macroscopic velocity, and $10^{-4}$ for the temperature, showing that the accuracy of RC does not change as $\Kn$ decreases. Moreover, the number of SI steps for RC and FSM is nearly the same for different $\Kn$, indicating that the approximation of RC model will not affect the convergence of the original Boltzmann equation. The CPU time for both RC and FSM listed in Tab. \ref{tab:cavity_time} illustrates that that of RC is much shorter than FSM for each $\Kn$, especially for the small Knudsen number $\Kn = 0.01$, where it is reduced by a factor of more than ten compared with that of FSM. For the 2D problems, the advantage of  RC is much more pronounced compared to that of the 1D problems. This may be due to the fact that the FSM collision operator has to be evaluated many more times in the 2D problems. Besides, the detailed CPU time for different Knudsen numbers for RC is listed in Tab. \ref{tab:cavity_time_RC}. It shows that for different $\Kn$, the cost of the SVD and constructing the collision model is of the same order, and the main difference for the CPU time is the source iteration, which is larger with smaller $\Kn$. This is also consistent with that of FSM.

\begin{table}[htbp]
\centering
{\footnotesize
\begin{threeparttable}
\renewcommand{\arraystretch}{1.15}
\begin{tabular*}{\textwidth}{@{\extracolsep{\fill}}r r c r c c c c@{}}
\toprule
\multirow{2}{*}{$\Kn$} & \multirow{2}{*}{method} & \multirow{2}{*}{steps}
& \multirow{2}{*}{CPU (s)} & \multicolumn{4}{c}{$L^2$ error} \\
\cmidrule(lr){5-8}
& & & & $\rho$ & $u_x$ & $u_y$ & $T$ \\
\midrule
\multirow{2}{*}{1}
& FSM & 10  & 14180 & \multirow{2}{*}{$3.74{\times}10^{-5}$} & \multirow{2}{*}{$1.51{\times}10^{-5}$} & \multirow{2}{*}{$1.64{\times}10^{-5}$} & \multirow{2}{*}{$1.62{\times}10^{-4}$} \\
& RC  & 10  & 3487   & & & & \\
\addlinespace
\multirow{2}{*}{0.1}
& FSM & 51  & 57873 & \multirow{2}{*}{$4.68{\times}10^{-5}$} & \multirow{2}{*}{$1.55{\times}10^{-5}$} & \multirow{2}{*}{$1.59{\times}10^{-5}$} & \multirow{2}{*}{$1.66{\times}10^{-4}$} \\
& RC  & 55  & 8427   & & & & \\
\addlinespace
\multirow{2}{*}{0.01}
& FSM & 684 & 823786 & \multirow{2}{*}{$3.31{\times}10^{-5}$} & \multirow{2}{*}{$5.66{\times}10^{-5}$} & \multirow{2}{*}{$3.12{\times}10^{-5}$} & \multirow{2}{*}{$1.69{\times}10^{-4}$} \\
& RC  & 662 & 74716  & & & & \\
\bottomrule
\end{tabular*}
\end{threeparttable}
}
\caption{(2D lid-driven cavity flow problem in Sec.~\ref{subsubsec:cavity}) 
The comparison of efficiency between RC and FSM for the 2D lid-driven cavity flow problem. Here, the CPU time for RC includes both the offline construction of the RC model and the online source iteration.}
\label{tab:cavity_time}
\end{table}

\begin{table}[!hptb]
\centering
{\footnotesize
\renewcommand{\arraystretch}{1.1}
\begin{tabular*}{0.75\textwidth}{@{\extracolsep{\fill}}r r c r r@{}}
\toprule
\multirow{2}{*}{$\Kn$} & \multicolumn{2}{c}{offline (s)} & \multirow{2}{*}{online SI (s)} & \multirow{2}{*}{total (s)} \\
\cmidrule(lr){2-3}
 & data generation & SVD and kernel & & \\
\midrule
1    & 512  & 2612 & 363   & 3487 \\
0.1  & 3872  & 2941 & 1614  & 8427 \\
0.01 & 47532 & 3820 & 23364 & 74716 \\
\bottomrule
\end{tabular*}
}
\caption{(2D lid-driven cavity flow problem in Sec.~\ref{subsubsec:cavity}) Detailed CPU time of the RC model in Tab.~\ref{tab:cavity_time} for different Knudsen numbers.}
\label{tab:cavity_time_RC}
\end{table}

\subsection{Time-dependent Boltzmann equation: fluid diffusion}
\label{sec:diffusion}
In this section, the time-dependent 2D fluid diffusion problem is studied. The computational domain is $\Omega=[0,1]\times[0,1]$, and the homogeneous Neumann boundary conditions are imposed on all sides to mimic an unbounded domain. The Knudsen number is fixed at $\Kn=0.03$. The initial condition is set as a Maxwellian 
\begin{equation}
    f(\bx,\boldsymbol{v},0)
    =
    \frac{\rho(\bx)}{(2\pi T)^{3/2}}
    \exp\!\left(-\frac{|\boldsymbol{v} - \bu|^2}{2T}\right),
\end{equation}
where the macroscopic variables are as follows
\begin{equation}
    \rho(\bx)=
    \begin{cases}
        10, & |x-0.5|\leqslant 0.05 \text{ and } |y-0.5|\leqslant 0.05,\\
        1, & \text{otherwise}.
    \end{cases}, \qquad \bu = \bm{0}, \qquad T = 1. 
\end{equation}
This initial condition generates a transient diffusion process with strong spatial gradients and a sharp central interface at early times. A similar fluid diffusion problem is also studied in \cite{cai2022numerical}.

The number of spatial cells is $N_x = N_y = 200$, where the MUSCL scheme \cite{vanleer1979towards} with the Minmod limiter \cite{sweby1984high} is utilized in the spatial discretization, and the time discretization is performed by a first-order semi-implicit scheme with the time step $\Delta t = 2.5\times 10^{-4}$. The number of Fourier modes is set as $N = 16$, and the rank number in RC is chosen as $N_f = N_Q = 120$ to account for the higher complexity of the transient dynamics with sharp gradients, and the reduced bases are constructed offline from the data generated on the coarse grid $N_x=N_y=20$.

\begin{figure}[!hptb]
\centering
\subfigure[$\rho$, $t=0.04$]{
            \includegraphics[width = 0.3\textwidth]{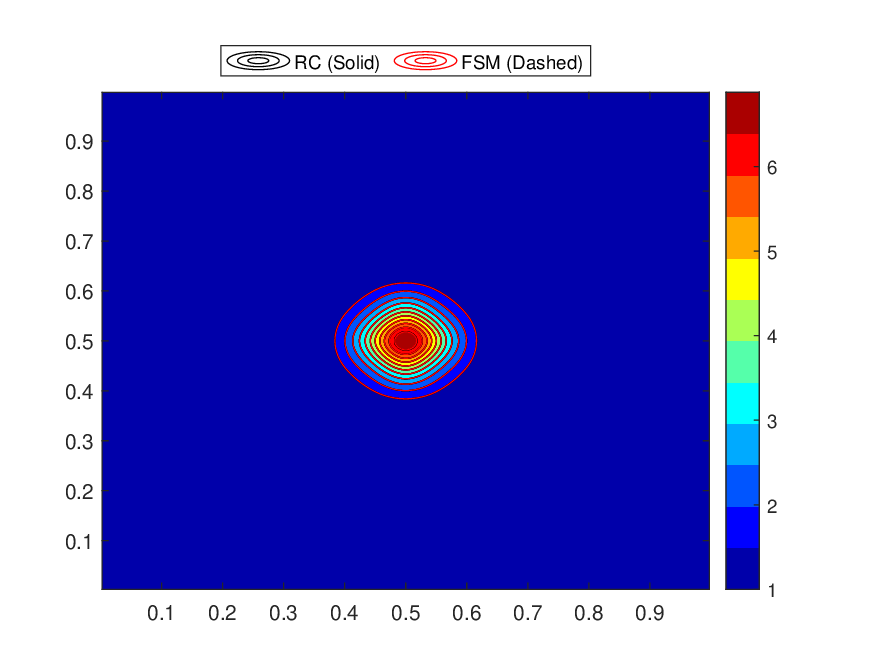}
}
\subfigure[$T$, $t=0.04$]{
            \includegraphics[width = 0.3\textwidth]{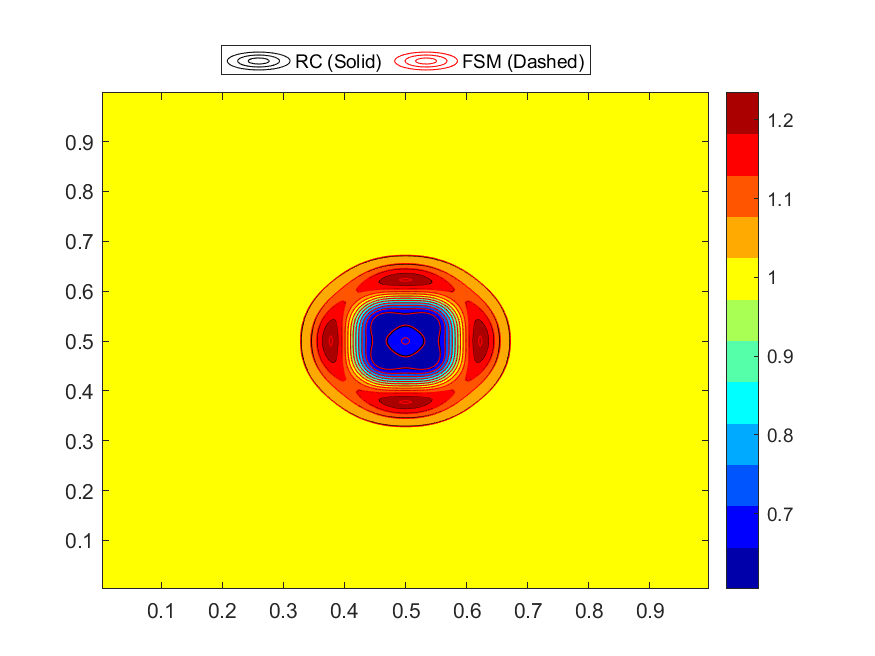}
}
\subfigure[$\sigma_{12}$, $t=0.04$]{
            \includegraphics[width = 0.3\textwidth]{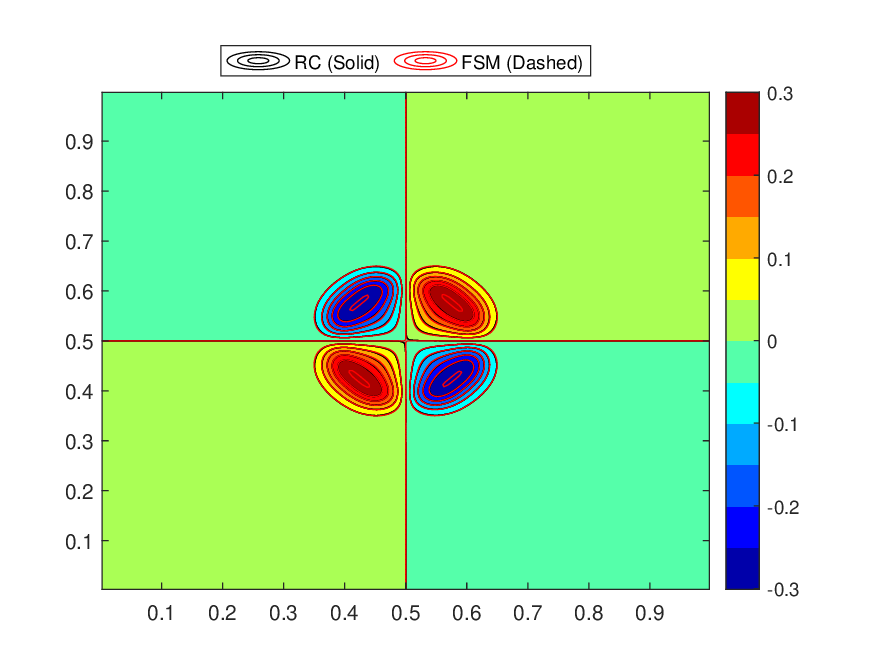}
} \\
\subfigure[$\rho$, $t=0.08$]{
            \includegraphics[width = 0.3\textwidth]{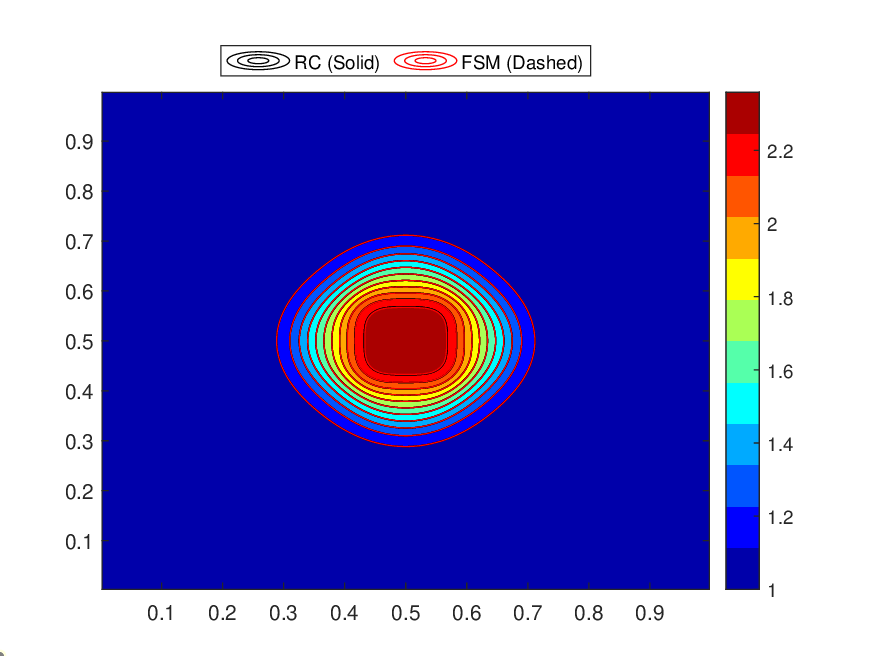}
}
\subfigure[$T$, $t=0.08$]{
            \includegraphics[width = 0.3\textwidth]{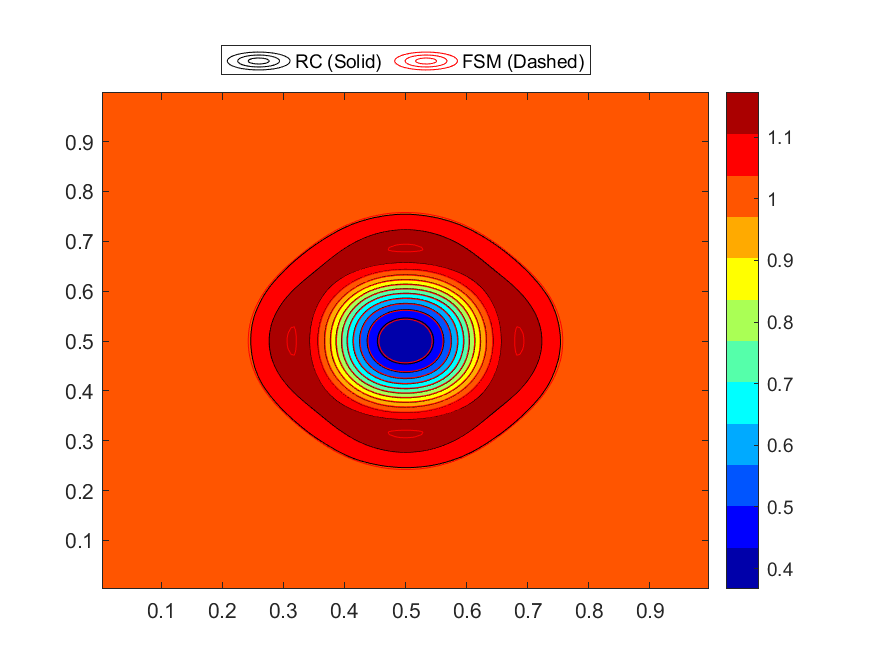}
}
\subfigure[$\sigma_{12}$, $t=0.08$]{
            \includegraphics[width = 0.3\textwidth]{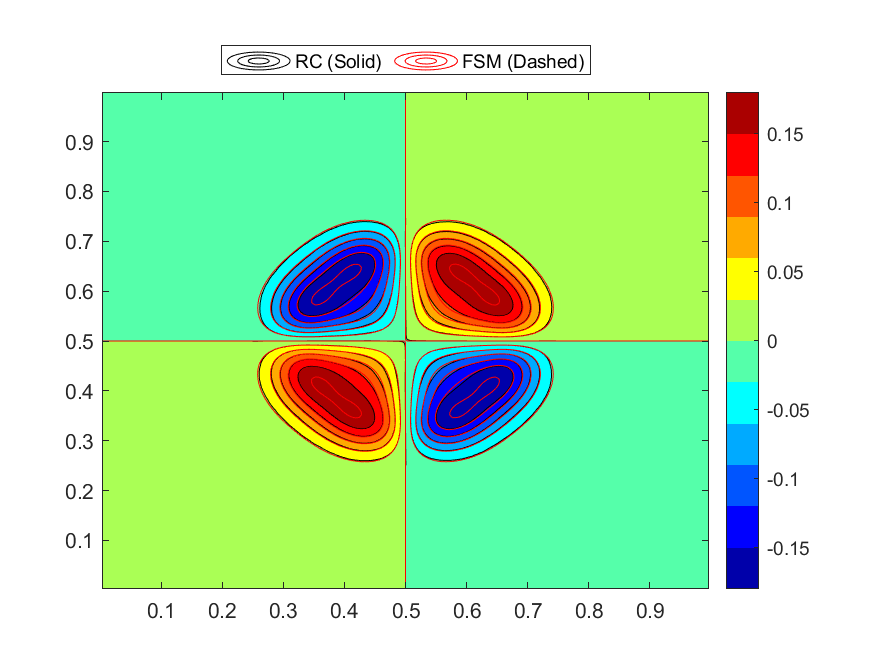}
} \\
\subfigure[$\rho$, $t=0.12$]{
            \includegraphics[width = 0.3\textwidth]{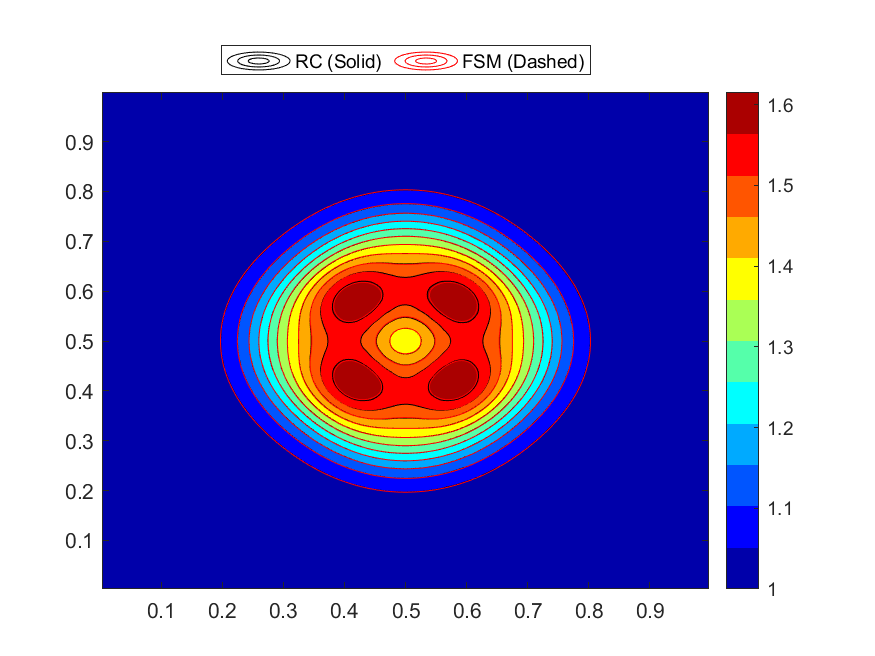}
}
\subfigure[$T$, $t=0.12$]{
            \includegraphics[width = 0.3\textwidth]{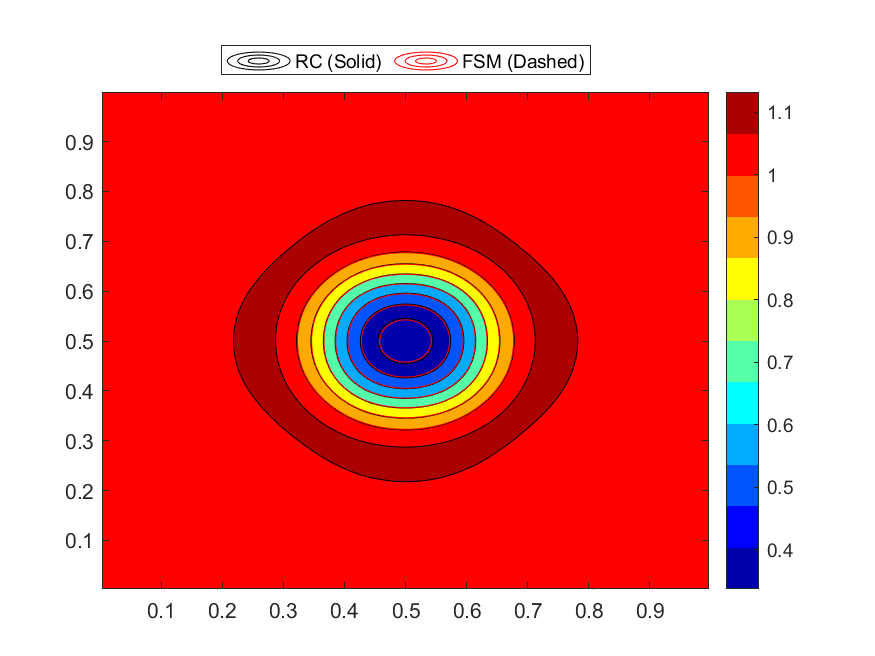}
}
\subfigure[$\sigma_{12}$, $t=0.12$]{
            \includegraphics[width = 0.3\textwidth]{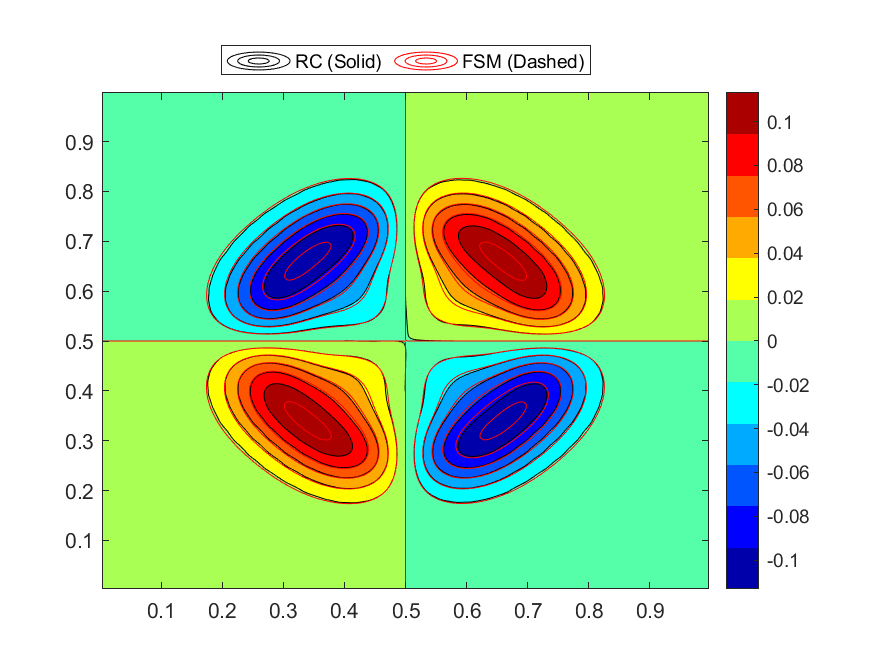}
}
\caption{(2D fluid-diffusion problem in Sec.~\ref{sec:diffusion}) Time evolution of the density $\rho$, the temperature $T$, and the shear stress $\sigma_{12}$ at $t=0.04,0.08$, and $0.12$. RC contours (solid black) are compared with FSM contours (red dashed).}
\label{fig:diffusion}
\end{figure}

Fig.~\ref{fig:diffusion} shows the numerical solution of the density $\rho$, the temperature $T$, and the shear stress $\sigma_{12}$ at $t = 0.04, 0.08$ and $0.12$. For all the cases, the numerical solution of RC matches well with that obtained by FSM. As time evolves, the high-density gas at the center diffuses into the surrounding region. The surrounding gas is compressed, and its temperature rises, while the temperature of the gas in the center decreases. The shear stress appears at four spots around the center, with opposite signs in neighboring quadrants, and it becomes weaker as the diffusion proceeds. This shows that the RC model remains stable and accurate in an unsteady problem with strong gradients.

Then, the computational efficiency of RC is examined. The $L^2$ error of the numerical solution between RC and FSM on the same mesh size $N_x = N_y = 200$ is calculated. This error of the density $\rho$ and the temperature $T$ at $t = 0.04, 0.08$ and $0.12$ are listed in Tab.~\ref{tab:diffusion_time}, where the GPU time of the offline and the online stages is also shown. It shows that the $L^2$ errors grow with the final time, since the error accumulates over the time steps. The GPU time for both RC and FSM listed in Tab.~\ref{tab:diffusion_time} illustrates that that of RC is much shorter than that of FSM for each final time, especially for $t = 0.12$, where it is reduced by a factor of more than eight compared with that of FSM. Moreover, the advantage of RC becomes more pronounced as the final time increases. This is because the offline time of RC, which is dominated by the SVD and the construction of the collision model, does not depend on the final time, while the online time grows linearly with the number of the time steps, so that the offline cost becomes less significant when more time steps are performed.

\begin{table}[htbp]
\centering
{\footnotesize
\renewcommand{\arraystretch}{1.15}
\begin{tabular*}{0.85\textwidth}{@{\extracolsep{\fill}}r r r r c c@{}}
\toprule
\multirow{2}{*}{$t$} & \multirow{2}{*}{method} & \multirow{2}{*}{offline (s)}
& \multirow{2}{*}{online (s)} & \multicolumn{2}{c}{$L^2$ error} \\
\cmidrule(lr){5-6}
& & & & $\rho$ & $T$ \\
\midrule
\multirow{2}{*}{0.04}
& FSM & --- & 6008 & \multirow{2}{*}{$3.99{\times}10^{-5}$} & \multirow{2}{*}{$2.85{\times}10^{-4}$} \\
& RC  & 571 & 553 & & \\
\addlinespace
\multirow{2}{*}{0.08}
& FSM & --- & 12016 & \multirow{2}{*}{$8.11{\times}10^{-5}$} & \multirow{2}{*}{$6.98{\times}10^{-4}$} \\
& RC  & 571 & 1107 & & \\
\addlinespace
\multirow{2}{*}{0.12}
& FSM & --- & 18024 & \multirow{2}{*}{$1.39{\times}10^{-4}$} & \multirow{2}{*}{$7.02{\times}10^{-4}$} \\
& RC  & 571 & 1661 & & \\
\bottomrule
\end{tabular*}
}
\caption{(2D fluid-diffusion problem in Sec.~\ref{sec:diffusion}) The comparison of efficiency between RC and FSM, where the GPU time of the RC model is split into the offline and the online stages.}
\label{tab:diffusion_time}
\end{table}

\section{Conclusion}
\label{sec:conclusion}
In this work, a reduced collision model is constructed for the Boltzmann equation to accelerate the evaluation of the quadratic collision model. For each target problem, the numerical solution on a coarser spatial mesh is utilized to construct the reduced basis for the distribution function and the collision model, where SVD is adopted, and techniques to preserve the conservation of the reduced collision are applied. Steady-state and time-dependent numerical examples are studied to validate this RC model. More numerical tests and the theoretical analysis of RC will be left for future work.

\section*{Funding}
The work of Yanli Wang was partially supported by the Science Challenge Program (No. TZ2026004). The work of Pei Zhang was supported by the China Scholarship Council (No. 202404890001).

\section*{Data Availability}
Enquiries about data availability should be directed to the authors.

\section*{Declarations}

\paragraph{Conflict of interest} The authors have not disclosed any competing interest.

\appendix
\section{Appendix}
\label{app:bilinear}
\subsection{Evaluation of the bilinear collision term by FSM}
In Sec. \ref{sec:dvm}, the evaluation of the quadratic collision term for the same discrete vector 
$\mQ(\boldsymbol{f},\boldsymbol{f})$ is listed in \eqref{eq:dis_hatQ}. Here, the calculation of $\mQ(\boldsymbol{g}_i,\boldsymbol{g}_j)$ in \eqref{eq:coe_Q} is presented, which also follows the fast Fourier spectral discretization in Sec.~\ref{sec:dvm}. 

For the continuous form, the bilinear collision operator is defined in the symmetric form as
\begin{equation}
\label{eq:app_bilinear}
    \mQ(a,b)=\frac{1}{2}\int_{\bbR^{3}}\int_{\mathbb{S}^{2}} B(\bv-\bv_{\ast},\sigma)\left[a(\bv_{\ast}^{\prime})b(\bv^{\prime})+a(\bv^{\prime})b(\bv_{\ast}^{\prime})-a(\bv_{\ast})b(\bv)-a(\bv)b(\bv_{\ast})\right]\dd\sigma\dd\bv_{\ast},
\end{equation}
where $a(\bv)$ and $b(\bv)$ are two continuous functions. Thus, $\mQ(a,b)=\mQ(b,a)$, and $\mQ(f,f)$ recovers the collision operator \eqref{eq:quic_coll}.

Let $\boldsymbol{a},\boldsymbol{b}\in\bbR^{N_v}$ be the discrete vectors of $a(\bv)$ and $b(\bv)$ on the discrete velocity points \eqref{eq:dis_v}, and  $\widehat{a}_{\boldsymbol{k}}$ and $\widehat{b}_{\boldsymbol{k}}$ be their Fourier coefficients calculated as in \eqref{eq:coe_f}. Similar to \eqref{eq:dis_hatQ}, the Fourier coefficients of the discrete bilinear collision term are evaluated as
\begin{equation}
\label{eq:app_hatQ}
\begin{aligned}
    \widehat{Q}_{\boldsymbol{k}}(\boldsymbol{a},\boldsymbol{b}) &=
    \frac{2\pi^2}{L^3}\sum_{n=1}^{J}\sum_{m=1}^{M}w_n^{\varrho} w_m^{\sigma} B(\varrho_n)
    \operatorname{sinc}\left(\frac{\pi \varrho_n |\boldsymbol{k}|}{2L}\right)I_{n,\boldsymbol{k},m}(\boldsymbol{a},\boldsymbol{b})   \\
    &\qquad -\frac{\pi^2}{L^3}\sum_{\boldsymbol{\ell}\in\mathcal{I}_N}
    \left(\widehat a_{\boldsymbol{k}-\boldsymbol{\ell}}\,\widehat b_{\boldsymbol{\ell}}+\widehat b_{\boldsymbol{k}-\boldsymbol{\ell}}\,\widehat a_{\boldsymbol{\ell}}\right)
    \left(\sum_{n=1}^{J} w_n^{\varrho} B(\varrho_n) \operatorname{sinc} \left(\frac{\pi \varrho_n |\boldsymbol{\ell}|}{L} \right) \right),
    \qquad \boldsymbol{k}\in\mathcal{I}_N,
\end{aligned}
\end{equation}
where
\begin{equation}
\label{eq:app_I}
\begin{aligned}
    I_{n,\boldsymbol{k},m}(\boldsymbol{a},\boldsymbol{b}) = \sum_{\boldsymbol{\ell}\in\mathcal{I}_N}
    \widehat a_{\boldsymbol{k}-\boldsymbol{\ell}}\,\widehat b_{\boldsymbol{\ell}}
    \bigg[ & \cos\left(\frac{\pi \varrho_n\boldsymbol{\sigma}_m
    \cdot (\boldsymbol{k}-\boldsymbol{\ell})}{2L} \right)
    \cos\left( \frac{\pi \varrho_n  \boldsymbol{\sigma}_m\cdot \boldsymbol{\ell}}{2L}  \right)  \\
    &+ \sin\left( \frac{\pi \varrho_n \boldsymbol{\sigma}_m\cdot
        (\boldsymbol{k}-\boldsymbol{\ell})}{2L} \right)
    \sin\left(  \frac{\pi \varrho_n \boldsymbol{\sigma}_m\cdot  \boldsymbol{\ell}}{2L}  \right) \bigg].
\end{aligned}
\end{equation}
It is easy to verify that $I_{n,\boldsymbol{k},m}(\boldsymbol{a},\boldsymbol{b})=I_{n,\boldsymbol{k},m}(\boldsymbol{b},\boldsymbol{a})$ by replacing $\boldsymbol{\ell}$ with $\boldsymbol{k}-\boldsymbol{\ell}$, and \eqref{eq:app_hatQ} reduces to \eqref{eq:dis_hatQ} when $\boldsymbol{a}=\boldsymbol{b}=\boldsymbol{f}$. Then, the discrete bilinear collision term is obtained by transforming $\widehat{Q}_{\boldsymbol{k}}(\boldsymbol{a},\boldsymbol{b})$ back to the velocity space as in \eqref{eq:dis_Q},
\begin{equation}
\label{eq:app_Q}
    Q_{\boldsymbol{p}}(\boldsymbol{a},\boldsymbol{b}) = \frac{1}{(2L)^3}\sum_{\boldsymbol{k}\in\mathcal{I}_N}
    \widehat Q_{\boldsymbol{k}}(\boldsymbol{a},\boldsymbol{b})\exp\left(\frac{\mathrm{i}\pi}{N}\boldsymbol{k}\cdot\boldsymbol{p}\right),
    \qquad \boldsymbol{p}\in\mathcal{I}_N.
\end{equation}
Finally, $\mQ(\boldsymbol{g}_i,\boldsymbol{g}_j)$ in \eqref{eq:coe_Q} is obtained by setting $\boldsymbol{a}=\boldsymbol{g}_i$ and $\boldsymbol{b}=\boldsymbol{g}_j$ in \eqref{eq:app_hatQ}--\eqref{eq:app_Q}, where the convolution sums are computed by the fast Fourier transform as in FSM \cite{gamba2017fast}. Due to the symmetry of \eqref{eq:app_bilinear}, the reduced kernel tensor \eqref{eq:kerner-tensor} satisfies $\mathcal{K}_{ijr}=\mathcal{K}_{jir}$.
\bibliographystyle{plain}
\bibliography{article}

\end{document}